\documentclass[dvipsnames,a4paper,twoside,reqno]{amsart}
\usepackage{amsthm,amsmath,amssymb,xcolor,xspace,bm}
\usepackage[bookmarks=true,bookmarksopen=true]{hyperref}
\theoremstyle{plain}
\newtheorem{theorem}{Theorem}[section]
\newtheorem{lemma}[theorem]{Lemma}
\newtheorem{corollary}[theorem]{Corollary}
\newtheorem{proposition}[theorem]{Proposition}

\theoremstyle{definition}

\newtheorem{assumption}[theorem]{Assumption}

\theoremstyle{remark}
\newtheorem{remark}[theorem]{Remark}

 \numberwithin{equation}{section}
 \numberwithin{table}{section}

\title[Discretization of bilinear control problems.]{Error estimates for the discretization of bilinear control problems governed by semilinear elliptic PDEs.}

\thanks{The first and third authors were supported by MCIN/ AEI/10.13039/501100011033/ under research project PID2020-114837GB-I00.}

\thanks{The second author was supported by the Hellenic Foundation for Research and Innovation (H.F.R.I.) under the ``First Call for H.F.R.I. Research Projects to support Faculty members and Researchers and the procurement of high-cost research equipment grant" (Project Number: 3270)}

\author[E. Casas]{Eduardo Casas}
\address{Departamento de Matem\'{a}tica Aplicada y Ciencias de la Computaci\'{o}n, E.T.S.I. Industriales y de Telecomunicaci\'on, Universidad de Cantabria, 39005 Santander, Spain.}
\email{eduardo.casas@unican.es}

\author[K. Chrysafinos]{Konstantinos Chrysafinos}
\address{Department of Mathematics, School of Applied Mathematics and Physical Sciences, National Technical University of Athens, Zografou Campus, Athens 15780, Greece and IACM, FORTH, 20013 Heraklion, Crete, Greece.}
\email{chrysafinos@math.ntua.gr}

\author[M. Mateos]{Mariano Mateos}
\address{Departamento de Matem\'{a}ticas, Campus de Gij\'on, Universidad de Oviedo, 33203, Gij\'on, Spain.}
 \email{mmateos@uniovi.es}



\ifpdf
\hypersetup{
  pdftitle={Discretization and error analysis of bilinear control problems governed by semilinear elliptic PDEs.},
  pdfauthor={E.~Casas,  K. Chrysafinos, and M.~Mateos}
}
\fi

\newcommand{\dx}{\,\mathrm{d}x}

\newcommand{\Pb}{\mbox{\rm (P)}\xspace}
\newcommand{\Pbh}{\mbox{\rm (P$_h$)}\xspace}

\newcommand{\uad}{U_{\rm ad}}
\newcommand{\Uh}{\mathcal{U}_h}
\newcommand{\uadh}{U_{h,\rm ad}}
\newcommand{\proj}{\operatorname{Proj}}

\newcommand{\A}{\mathcal{A}}

\newcommand{\dimension}{d}
\newcommand{\conormal}{n\!}
\newcommand{\tichonov}{\nu}

\newcommand{\umin}{{\alpha}}
\newcommand{\umax}{{\beta}}

\begin{document}

\begin{abstract}
This paper studies an optimal control problem governed by a semilinear elliptic equation, in which the control acts in a multiplicative or bilinear way as the reaction coefficient of the equation.
We focus on the numerical discretization of the problem.
The discretization is carried out by using the finite element method, with piecewise constant functions for the control and continuous piecewise linear functions for the state and the adjoint state. We first prove convergence of the solutions of the discrete problems to solutions of the continuous problem. We also demonstrate that strict local solutions of the continuous problem can be approximated by local solutions of the discrete problems. Next we obtain an error estimate of order $O(h)$ for the difference between continuous and discrete locally optimal controls. To obtain this result we assume no-gap second order sufficient optimality conditions. As it is usual in this kind of discretization, a superconvergence phenomenon of order $O(h^2)$ is observed in numerical experiments for the error estimates of the state and adjoint state. The last part of the paper is dedicated to explain this behaviour. A numerical experiment confirming these results is included.
\end{abstract}

\maketitle

\section{Introduction} \label{S1}
We study the finite element discretization of the following bilinear optimal control problem:
\[ \Pb \min_{u \in \uad} J(u) :=  \int_\Omega L(x,y_u(x)) \dx + \frac{\tichonov}{2} \int_\Omega u^2(x) \dx, \]
where $y_u$ is the state associated with the control $u$ solution of

\begin{equation}
\left\{\begin{array}{l} Ay + a(x,y) + uy= 0\ \  \mbox{in } \Omega,\vspace{2mm}\\  \partial_{\conormal_A} y = g\ \ \mbox{on }\Gamma. \end{array}\right.
\label{E1.1}
\end{equation}
Here $\Omega \subset \mathbb{R}^{\dimension}$, $\dimension =2$ or $3$, is a bounded open convex set with a  polygonal or polyhedral boundary $\Gamma$, $A$ is an elliptic operator, $L:\Omega\times \mathbb{R}\to\mathbb{R}$, $a:\Omega\times \mathbb{R}\to\mathbb{R}$ and $g:\Gamma\to\mathbb{R}$ are given functions, $\tichonov >0$ and the admissible set of controls is defined as
\[ \uad = \{ u \in L^2(\Omega) :  \umin \leq u(x) \leq \umax \text{ a.e. in }  \Omega \}, \]
where $-\infty < \umin<\umax\leq \infty$. The precise assumptions on the data $A$, $L$, $a$, $g$ and $\umin$ will be given in Sections \ref{S2} and \ref{S3}.

Kr\"oner and Vexler \cite{KronerVexler2009} and Winkler \cite{Winkler2020} have studied the finite element discretization of \Pb in the case of having a linear elliptic equation. In the second reference, the control acts as the Robin coefficient on the boundary. In bilinear control problems, even in case of linear equations, the control-to-state mapping is non-linear, so a second order sufficient condition is used to obtain error estimates. In both references, the authors assume coercivity of the second derivative of the cost functional at the solution with respect to all directions of the control space, which implies the strict convexity of the control problem in a neighborhood of the optimal control. We aim to obtain similar results removing this convexity assumption and only assuming no-gap second order sufficient optimality conditions.  Under an usual hypothesis, see Assumption \ref{WWA5.1}, we obtain a superconvergence result for the optimal state and adjoint state variables. As a consequence, we also get order $O(h^2)$ for the approximation of the optimal control by using the post-processing technique proposed by Meyer and R\"osch \cite{Meyer-Rosch2004}. The reader is also referred to \cite{CWW2018} for the proof of error estimates in the case where the Tikhonov term is not present ($\nu = 0$). Finally, let us mention that we also allow non bounded controls from above and a small negative value for the lower bound.

The plan of the paper is as follows. In the next section, we summarize the results obtained in \cite{Casas-Chrysafinos-Mateos2023} about the continuous problem and establish some regularity results that are needed to prove error estimates. In Section \ref{S3}, we discretize the control problem and prove the convergence of the discretizations following the results of \cite{Casas-Troltz2012COAP}. Section \ref{WWS4} is devoted to the derivation of error estimates for the control. We discretize the control using piecewise constant functions and, consequently, convergence of order $O(h)$ is obtained. As it is shown in the numerical experiments in Section \ref{WWS6}, a superconvergence phenomenon in the state and adjoint state variable appears. This is completely justified in Section \ref{WWS5}, where order $O(h^2)$ is obtained for both variables under a natural assumption.

\section{Analysis of the continuous problem}\label{S2}
In our recent work \cite{Casas-Chrysafinos-Mateos2023}, we studied the continuous problem \Pb. We summarize below the results needed for the subsequent numerical analysis. For the proofs, the reader is referred to the above mentioned paper. The following hypotheses are assumed along this paper.
\begin{assumption} \label{A2.1}
The operator $A$ defined in $\Omega$ and the conormal derivative $\partial_{n_A}$ on $\Gamma$ are given by the expressions
\[
Ay = - \sum_{i,j=1}^{\dimension} \partial_{x_j} [ a_{ij}(x) \partial_{x_i} y] \quad \text{ and } \quad  \partial_{n_A}y = \sum_{i,j = 1}^da_{ij}(x)\partial_{x_i}y(x)n_j(x),
\]
where $n(x)$ denotes the outward unit normal vector to $\Gamma$ at the point $x$ and $a_{ij} \in C^{0,1}(\bar\Omega)$ for $1\leq i,j \leq \dimension$ satisfies for some $M_A, \Lambda_A>0$
\[ M_A | \xi |^2 \geq \sum_{i,j=1}^{\dimension} a_{ij}(x) \xi_i \xi_j \geq \Lambda_A | \xi |^2 \quad \forall x \in \Omega \ \text{and} \ \forall \xi \in \mathbb R^n.
\]
\end{assumption}

\begin{assumption} \label{A2.2}
We assume that $a: \Omega \times \mathbb{R} \longrightarrow \mathbb{R}$ is a Carath\'eodory function of class $C^2$ with respect to the second variable satisfying the following properties for a.a. $x \in \Omega$:
\begin{align*}
&\bullet a(\cdot,0) \in L^{2}(\Omega), \\
&\bullet  \exists a_0 \in L^\infty(\Omega) \text{ such that }\frac{\partial a}{\partial y}(x,y) \geq a_0(x) \ge 0 \ \forall y \in \mathbb{R}, \\
&\bullet \forall M > 0 \ \exists C_{a,M} \text{ such that } \sum_{j=1}^2 \Big| \frac{\partial^j a}{\partial y^j} (x,y) \Big| \leq C_{a,M}  \  \forall | y | \leq M, \\
&\bullet \forall \varepsilon >0 \text{ and } \forall M>0 \ \exists \rho >0 \text{ such that } \ \Big|    \frac{\partial^2 a}{\partial y^2}  (x,y_1) -  \frac{\partial^2 a}{\partial y^2}  (x,y_2) \Big| \leq \varepsilon \\
& \qquad \qquad \text{ for all } |y_1|, |y_2| \leq M \text{ with } |y_1-y_2 | \leq \rho.
\end{align*}
\end{assumption}

\begin{assumption} \label{A2.3}
The boundary data satisfies the following regularity property: $g=\partial_{\conormal_A} y_g$ for some $y_g \in H^2(\Omega)$.
\end{assumption}

\begin{remark} \label{WWR2.4} We observe that the main results of \cite{Casas-Chrysafinos-Mateos2023}  were proven under the less stringent assumptions of $a_{ij} \in L^\infty(\Omega)$, $a(\cdot,0) \in L^p(\Omega)$ for some $p > \frac{d}{2},$ and $g \in L^{q}(\Gamma)$ for some $q > d-1.$ The enhanced regularity assumptions are necessary to establish $H^2(\Omega)$ regularity for the solution of \eqref{E1.1}.
In our case, since $y_g\in H^2(\Omega)$ and the coefficients $a_{ij}$ are Lipschitz, from the trace theorem and Sobolev's embedding theorem we have that $g\in L^4(\Gamma)$ for $d\leq 3$.
\end{remark}

Throughout this paper the following notation will be used:
\[
 m_u := \text{ess}\,  \text{inf}_{x \in \Omega} u(x), \ \A_0 := \{ u \in L^2(\Omega) : a_0(x) + m_u \ge 0 \text{ a.e. in } \Omega \text{ and } a_0+u\not\equiv 0\}.
\]
From Assumption \ref{A2.1}, for every $u\in\A_0$ we infer the existence of a constant $0 < \Lambda_u \le \Lambda_A$ such that
\begin{equation}
\int_\Omega\Big(\sum_{i,j=1}^{\dimension}a_{i,j}(x)\partial_{x_i} y(x)\partial_{x_j}y(x) + [a_0(x) + u(x)]y(x)^2\Big)\dx \ge \Lambda_u\Vert y\Vert ^2_{H^1(\Omega)} \ \ \forall y \in H^1(\Omega).
\label{E2.1}
\end{equation}
The following theorem summarizes existence, uniqueness and regularity results for \eqref{E1.1}.
\begin{theorem} \label{WWT2.5}
For every $u \in \A_0$ there exists a unique solution $y_u \in H^2(\Omega)$ of \eqref{E1.1}. Furthermore, the following estimates hold:
\begin{align}
& \Vert y_u\Vert_{H^1(\Omega)} \leq \frac{1}{\Lambda_u} \left ( C_{\Omega} \Vert a(\cdot,0)\Vert_{L^2(\Omega)} + C_{\Gamma} \Vert g\Vert_{L^4(\Omega)} \right ), \label{E2.2} \\
& \Vert y_u\Vert_{L^{\infty}(\Omega)} \leq M_\infty , \label{E2.3}
\\
&\Vert y_u\Vert_{H^2(\Omega)} \leq C_2 \left ( 1 + \Vert u\Vert_{L^2(\Omega)}  \right ). \label{E2.4}
\end{align}
where $M_{\infty}$ depends on $\Lambda_u$, $\Vert a(\cdot,0)\Vert_{L^2(\Omega)}$ and  $\Vert g\Vert_{L^4(\Gamma)}$, and $C_{2}$ depends on $M_\infty$ and $\Vert y_g\Vert_{H^2(\Omega)}$.

\end{theorem}
\begin{proof}
Existence and uniqueness in $H^1(\Omega) \cap L^\infty(\Omega)$ as well as estimates \eqref{E2.2}, \eqref{E2.3} are obtained in \cite[Theorem 2.4]{Casas-Chrysafinos-Mateos2023}.
  Let us prove the $H^2(\Omega)$ regularity of $y_u$. Denote $F(x) = (1- u(x))y_u(x) - a(x,y_u(x))$. From Assumptions \eqref{A2.2} and \eqref{E2.3}, we deduce that $F\in L^2(\Omega)$ and
\[
\Vert F\Vert_{L^2(\Omega)} \leq (|\Omega|^{\frac{1}{2}} + \Vert u\Vert_{L^2(\Omega)})M_\infty + \Vert a(\cdot,0)\Vert_{L^2(\Omega)}+ C_{a,M_\infty} M_\infty|\Omega|^{\frac{1}{2}}.
\]
  The function $z = y_u-y_g$ satisfies
  \[
  \left\{\begin{array}{l} Az + z= F - Ay_g - y_g\ \  \mbox{in } \Omega,\vspace{2mm}\\  \partial_{\conormal_A} z = 0\ \ \mbox{on }\Gamma. \end{array}\right.
\]
From Assumptions \ref{A2.1} and \ref{A2.3}, we readily deduce that $Ay_g+y_g \in L^2(\Omega)$ and the existence of some constant $C>0$ such that
\[
\Vert Ay_g+y_g\Vert_{L^2(\Omega)}\leq C\Vert y_g\Vert_{H^2(\Omega)}.
\]
The $H^2(\Omega)$ regularity of $z$ is a consequence of the Lipschitz regularity of the coefficients $a_{ij}$ and the convexity of $\Omega$; see \cite[Chapter 3]{Grisvard85}. Moreover, we also have the estimate
\[
\Vert z\Vert_{H^2(\Omega)} \le C'\Big(\Vert Ay_g+y_g\Vert_{L^2(\Omega)} + \Vert F\Vert_{L^2(\Omega)}\Big)
\]
for a constant $C'$ depending only on $A$. The $H^2(\Omega)$ regularity of $y_u$ follows from this and the regularity of $y_g$ stated in Assumption \ref{A2.3}. The estimate \eqref{E2.4} follows from the estimates for $z$ and $F$.
\end{proof}

The following result is of utmost importance. It establishes the differentiability of the control-to-state mapping in an open subset  $\A\subset L^2(\Omega)$ containing $\A_0$.
\begin{theorem} \label{WWT2.6}
There exists an open set $\A$ in $L^2(\Omega)$ such that $\A_0 \subset \A$ and $\forall u \in \A$ the equation \eqref{E1.1} has a unique solution $y_u \in H^2(\Omega)$.
Further, the mapping $G : \A \longrightarrow H^2(\Omega)$ defined by $G(u):=y_u$ is of class $C^2$ and $\forall u \in \A$ and
$\forall v \in L^2(\Omega)$ the function $z_{u,v} =G'(u)v$ is the unique solution of the equation:
\begin{align}
& \left\{\begin{array}{l} Az + \displaystyle \frac{\partial a}{\partial y}(x,y_u)z + uz =- vy_u\ \  \mbox{in } \Omega,\vspace{2mm}\\  \partial_{\conormal_A} z = 0\ \ \mbox{on }\Gamma. \end{array}\right.
\label{E2.5}
\end{align}
In addition, $\forall u \in \A$ and $\forall v_1, v_2 \in L^2(\Omega)$, the function $w=G''(u)(v_1,v_2)$ is the unique solution of the equation:
\begin{align*}
& \left\{\begin{array}{l} Aw + \displaystyle \frac{\partial a}{\partial y}(x,y_u)w + uw = -\Big[\frac{\partial^2a}{\partial y^2}(x,y_u)z_{v_1}z_{v_2} + v_1z_{u,v_2} + v_2z_{u,v_1}\Big] \ \  \mbox{in } \Omega,\vspace{2mm}\\  \partial_{\conormal_A} w = 0\ \ \mbox{on }\Gamma, \end{array}\right.
\nonumber
\end{align*}
where $z_{u,v_i} = G'(u)v_i$, $i=1,2$.
\end{theorem}
\begin{proof}
Let us consider the space $Y_\Gamma = \{\partial_{n_A}y : y \in H^2(\Omega)\}$ endowed with the norm $\Vert h\Vert_{Y_\Gamma} = \inf\{\Vert y\Vert_{H^2(\Omega)} : \partial_{n_A}y = h\}$. Then, $Y_\Gamma$ is a Banach space. We define the function
\[
\mathcal{F}:L^2(\Omega) \times H^2(\Omega) \longrightarrow L^2(\Omega) \times Y_\Gamma,\ \ \mathcal{F}(u,y) = (Ay + a(x,y) + uy,\partial_{n_A}y - g).
\]
Following the steps of \cite[Theorem 2.5]{Casas-Chrysafinos-Mateos2023}, the proof of the theorem is achieved by applying the implicit function theorem to the function $\mathcal{F}$.
\end{proof}
Next we recall the main assumptions and results about the control problem \Pb already established in \cite{Casas-Chrysafinos-Mateos2023}.
\begin{assumption} \label{WWA2.7} We assume that $a_0(x) + \umin \ge 0$ for a.a. $x \in \Omega$, $a_0 + \umin \not\equiv 0$, and $\umin < \umax \leq \infty$ hold.
\end{assumption}
\begin{assumption} \label{WWA2.8} The function $L: \Omega \times \mathbb R \longrightarrow \mathbb R$ is Carath\'eodory and of class of $C^2$ with respect to the second variable.
Further the following properties hold for almost all $x \in \Omega$:
\begin{align*}
& \bullet L(\cdot, 0) \in L^1(\Omega) \\
& \bullet \forall M >0, \  \exists L_M \in L^2(\Omega) \text{ such that } \Big|  \frac{\partial L}{\partial y}(x,y) \Big| \leq L_M(x) \ \forall |y| \leq M, \\
& \bullet \forall M >0, \ \exists C_{L,M} \in \mathbb R   \text{ such that } \Big|  \frac{\partial^2 L}{\partial y^2}(x,y) \Big| \leq C_{L,M} \ \forall |y| \leq M, \\
& \bullet \forall \varepsilon >0  \text{ and } \forall M >0 \ \exists \rho >0 \text{ such that } \\
& \qquad\qquad  \Big|  \frac{\partial^2 L}{\partial y^2}(x,y_1) - \frac{\partial^2 L}{\partial y^2}(x,y_2)  \Big| \leq \varepsilon \ \forall |y_1|, |y_2| \leq M \text{ with } |y_1-y_2| \leq \rho.
\end{align*}
\end{assumption}

As a consequence of Theorem \ref{WWT2.6} and Assumption \ref{WWA2.8} we deduce the differentiability of functional $J$.
\begin{theorem} \label{WWT2.9} The functional $J: \A \longrightarrow \mathbb R$ is of class $C^2$ and its derivatives are given by the expressions:
\begin{align}
 J'(u)v =& \int_\Omega (\tichonov u - y_u \varphi_u)v \dx \ \ \forall u \in \A , \ \forall v \in L^2(\Omega), \label{WWE2.6} \\
 J''(u)(v_1,v_2) =& \int_\Omega \Big[ \frac{\partial^2 L}{\partial y^2}(x,y_u) - \varphi_u  \frac{\partial^2 a}{\partial y^2}(x,y_u) \Big] z_{u,v_1} z_{u,v_2}\dx, \label{WWE2.7}\\
 & - \int_\Omega \Big[ v_1z_{u,v_2} + v_2 z_{u,v_1} \Big] \varphi_u \dx + \tichonov \int_\Omega v_1 v_2 \dx \nonumber\\
 = & \int_\Omega \left[\tichonov v_1 -(\varphi_u z_{u,v_1}+y_u\eta_{u,v_1})\right]v_2\dx\nonumber\\
 = & \int_\Omega \left[\tichonov v_2 -(\varphi_u z_{u,v_2}+y_u\eta_{u,v_2})\right]v_1\dx , \ \forall u \in \A , \ \forall v_1, v_2 \in L^2(\Omega), \nonumber
\end{align}
where $z_{u,v_i} = G'(u)v_i$, $i=1,2$,  $\varphi_u \in H^2(\Omega)$ is the adjoint state, the unique solution of the equation
\begin{equation}
 \left\{\begin{array}{l} A^*\varphi + \displaystyle \frac{\partial a}{\partial y}(x,y_u)\varphi + u\varphi = \frac{\partial L}{\partial y}(x,y_u) \ \  \mbox{in } \Omega,\vspace{2mm}\\  \partial_{\conormal_{A^*}} \varphi = 0\ \ \mbox{on }\Gamma, \end{array}\right.
\label{WWE2.8}
\end{equation}
and $\eta_{u,v_i}\in H^2(\Omega)$, $i = 1, 2$, denotes the solution of
\begin{equation}
 \left\{\begin{array}{l} A^*\eta + \displaystyle \frac{\partial a}{\partial y}(x,y_u)\eta + u\eta = \left[\frac{\partial^2 L}{\partial y^2}(x,y_u)-\varphi_u \frac{\partial^2 a}{\partial y^2}(x,y_u)\right] z_{u,v_i}-v_i\varphi_u\ \  \mbox{in } \Omega,\vspace{2mm}\\  \partial_{\conormal_{A^*}} \eta = 0\ \ \mbox{on }\Gamma, \end{array}\right.
\label{WWE2.9}
\end{equation}
\end{theorem}

\begin{remark} \label{WWR2.10}
1. From Assumption \ref{WWA2.7} and the fact that that $a_0(x) + m_u \geq a_0(x) + \umin \ \forall u \in \uad$ we infer that $\uad \subset \A_0 \subset \A$. \vspace{2mm}\\
2. Using again Assumption \ref{WWA2.7} we deduce the existence of a constant $\Lambda > 0$ such that the inequality \eqref{E2.1} holds with $u$ and $\Lambda_u$ replaced by $\alpha$ and $\Lambda$. Since $a_0(x) + m_u \geq a_0(x) + \umin \ \forall u \in \uad$, without loss of generality we can assume that $\Lambda_u \ge \Lambda$ $\forall u \in \uad$. Consequently, \eqref{E2.2} and \eqref{E2.3} hold with $\Lambda$ instead of $\Lambda_u$, and $C_2$ in \eqref{E2.4} can be chosen independently of $u \in \uad$.
Arguing similarly we deduce the existence of $C_\infty$ such that
\[
\Vert \varphi_u\Vert_{L^\infty(\Omega)} + \Vert \varphi_u\Vert_{H^1(\Omega)} \le C_\infty \quad \forall u \in \uad.
\]
 \end{remark}

Now, we address the existence of a solution of \Pb and the first order optimality conditions. To this end, we will use the following lemma.
\begin{lemma}\label{WWL2.11}
Suppose that $\{u_k\}_{k=1}^\infty\subset \uad$. The following statements are fulfilled:

1- If $u_k\rightharpoonup u$ weakly in $L^2(\Omega)$ then the convergences $y_{u_k}\rightharpoonup y_u$ in $H^2(\Omega)$ and $y_{u_k} \to y_u$ in $C(\bar\Omega)$ hold.

2- If $\Vert u_k\Vert_{L^2(\Omega)} \to \infty$ then $J(u_k) \to \infty$ is satisfied.
\end{lemma}
\begin{proof}
1- Since $\uad$ is closed and convex in $L^2(\Omega)$, it is weakly closed and hence $u\in\uad$. Denote $\tilde y_k = y_u - y_{u_k}$. Using the mean value theorem we infer the existence of a measurable function $\theta_k:\Omega \longrightarrow [0,1]$ such that for $y_{\theta_k} = y_u + \theta_k(y_{u_k} - y_u)$ we have that
\begin{equation}\label{WWE2.10}
A\tilde y_k+\frac{\partial a}{\partial y}(x,y_{\theta_k}) \tilde y_k + u_k \tilde y_k = (u_k-u) y_{u}\mbox{ in }\Omega,\qquad \partial_{n_A} \tilde y_k=0.
\end{equation}
Testing this equation with $\tilde y_k$ and using Remark \ref{WWR2.10} and the compactness of the embedding $L^2(\Omega) \hookrightarrow H^{-1}(\Omega)$ we infer
\[
\Lambda\Vert \tilde y_k\Vert ^2_{H^1(\Omega)} \le \Vert y_u\tilde y_k\Vert_{H^1(\Omega)}\Vert u_k - u\Vert_{H^{-1}(\Omega)} \le C\Vert y_u\Vert_{H^2(\Omega)}\Vert u_k - u\Vert_{H^{-1}(\Omega)}\Vert \tilde y_k\Vert_{H^1(\Omega)}.
\]
This leads $\Vert \tilde y_k\Vert_{H^1(\Omega)} \le \frac{C}{\Lambda}\Vert y_u\Vert_{H^2(\Omega)}\Vert u_k - u\Vert_{H^{-1}(\Omega)} \to 0$ as $k \to \infty$. Further, from \eqref{E2.3} and Remark \ref{WWR2.10}-2, we have that $\Vert y_{u_k}\Vert_{L^\infty(\Omega)}\leq M_{\infty}$ for every $k$. The same estimate holds for $y_u$. Then, from Assumption \ref{A2.2} and the above equation we get that $A\tilde y_k$ is bounded in $L^2(\Omega)$ and, hence, $\{\tilde y_k\}_{k = 1}^\infty$ is bounded in $H^2(\Omega)$. Consequently, we have that $\tilde y_k \rightharpoonup 0$ in $H^2(\Omega)$. Finally, the compactness of the embedding $H^2(\Omega) \hookrightarrow C(\bar\Omega)$ implies the strong convergence $\tilde y_k \to 0$ in $C(\bar\Omega)$, which concludes the proof of 1.

2- According to \eqref{E2.3} and Remark \ref{WWR2.10}, we have that $\Vert y_{u_k}\Vert_{L^\infty(\Omega)} \le M_\infty$ for every $k \ge 1$. Now, from Assumption \ref{WWA2.8} we infer with the mean value theorem for all $k \ge 1$
\[
|L(x,y_{u_k}(x))| \le |L(x,0)| + L_{M_\infty}(x)M_\infty \Rightarrow \exists C > 0 \text{ such that} \int_\Omega |L(x,y_{u_k}(x))|\dx \le C.
\]
It is enough to observe that $-C + \frac{\nu}{2}\Vert u_k\Vert ^2_{L^2(\Omega)} \le J(u_k)$ to deduce that $J(u_k) \to \infty$.
\end{proof}

We will say that $\bar u\in\uad$ is a local minimum of $\Pb$ if there exists $\rho >0$ such that $J(\bar u)\leq J(u)$ for all $u\in\uad\cap B_\rho(\bar u)$, where $B_\rho(\bar u) = \{u\in L^2(\Omega):\Vert u-\bar u\Vert_{L^2(\Omega)}<\rho\}$. It is a strict local minimum if the inequality is strict for $u\neq \bar u$. In the sequel, $\bar B_\rho(\bar u)$ denotes the closed ball of $B_\rho(\bar u)$.

\begin{theorem} \label{WWT2.12} Problem $\Pb$ has at least one solution. Moreover, if $\bar u \in \uad$ is a local minimizer of $\Pb$ then there exist $\bar y, \bar \varphi \in H^2(\Omega)$ such that
\begin{align}
& \left\{\begin{array}{l} A\bar y + a(x,\bar y) + \bar u \bar y= 0\ \  \mbox{in } \Omega,\vspace{2mm}\\  \partial_{\conormal_A} \bar y = g\ \ \mbox{on }\Gamma. \end{array}\right.
\label{WWE2.11} \\
& \left\{\begin{array}{l} A^* \bar \varphi + \displaystyle \frac{\partial a}{\partial y}(x,\bar y)\bar \varphi + \bar u \bar \varphi = \frac{\partial L}{\partial y}(x,\bar y) \ \  \mbox{in } \Omega,\vspace{2mm}\\  \partial_{\conormal_{A^*}} \bar  \varphi = 0\ \ \mbox{on }\Gamma. \end{array}\right. \label{WWE2.12} \\
& \bar u(x) = \proj_{[\umin, \umax]} \left ( \frac{1}{\tichonov} \bar y(x) \bar \varphi(x) \right ). \label{WWE2.13}
\end{align}
Moreover, $\bar u\in  H^1(\Omega) \cap C^{0,\mu} (\bar \Omega)$ for all $\mu \in (0,1)$ if $\dimension =2$ and $\mu=1/2$ if $\dimension = 3$.
\end{theorem}

\begin{proof}
  Existence of a solution follows in standard way using the direct method of the calculus of variations with the help of Lemma \ref{WWL2.11}. First order optimality conditions are a consequence of the differentiability of $J$ in $\A$ and \eqref{WWE2.6}. Finally, we notice that the product of two functions in $H^2(\Omega)$ is also in $H^2(\Omega)$. Therefore, from the Sobolev embedding theorem we have that $\bar y\bar\varphi \in H^2(\Omega) \subset H^1(\Omega) \cap C^{0,\mu} (\bar \Omega)$ for all $0<\mu < 1$ if $\dimension =2$ and $\mu=1/2$ if $\dimension = 3$. Hence, the regularity of $\bar u$ is a straightforward consequence of \eqref{WWE2.13}.
\end{proof}

From now on $(\bar u,\bar y,\bar \varphi) \in \uad \times [ H^2(\Omega)]^2$ will denote a triplet that satisfies \eqref{WWE2.11}-\eqref{WWE2.13}.  Associated to this triplet we define the standard cone of critical directions,
\begin{equation} \label{WWE2.14}
\mathcal{C}_{\bar u} = \{ v \in L^2(\Omega) : v(x) {=} 0 \text{ if } \tichonov \bar u(x) - \bar y(x) \bar \varphi(x) \neq 0  \text{ a.e. in $\Omega$ and (\ref{WWE2.15})  holds} \},
\end{equation}
where
\begin{equation} \label{WWE2.15}
v(x)  \left\{ \begin{array}{ll} \geq 0 & \text{ if } \bar u(x) = \umin, \\ \leq 0 & \text{ if } \bar u(x) = \umax. \end{array} \right.
\end{equation}
In addition, we define for $\tau\geq 0$ the following extended cone of critical directions

\begin{equation*}
{\mathcal C}^{\tau}_{\bar u} = \{ v \in L^2(\Omega) : v(x) {=} 0 \text{ if } |\tichonov \bar u(x) - \bar y(x) \bar \varphi(x)| > \tau  \text{ a.e. in $\Omega$ and (\ref{WWE2.15})  holds} \}.
\end{equation*}

 Regarding the second order optimality conditions we have the following result.
\begin{theorem} \label{WWT2.15}
If $\bar u$ is a local minimizer of $\Pb$, then $J''(\bar u)v^2 \geq 0 \ \forall v \in \mathcal{C}_{\bar u}$ holds. Conversely, if $\bar u \in \uad$ satisfies the first order optimality conditions \eqref{WWE2.11}--\eqref{WWE2.13} and
\begin{equation}\label{WWE2.16}J''(\bar u)v^2 > 0 \ \forall v \in \mathcal{C}_{\bar u} \setminus\{0\},\end{equation}
then there exist $\varepsilon >0$ and $\delta >0$ such that
\begin{equation} \label{WWE2.17}
J(\bar u) + \frac{\delta}{2} \Vert u - \bar u\Vert ^2_{L^2(\Omega)} \leq J(u) \quad\forall u \in \uad \text{ with } \Vert u-\bar u\Vert_{L^2(\Omega)} \leq \varepsilon.
\end{equation}
\end{theorem}
This is a classical result whose proof can be found, for instance, in \cite{Casas-Troltzsch12}. The next theorem establishes an important equivalence that will be used later in the numerical analysis.

\begin{theorem}
Let $\bar u \in \uad$ satisfy the first order optimality conditions \eqref{WWE2.11}--\eqref{WWE2.13}. Then, $J''(\bar u)v^2 > 0 \ \forall v \in \mathcal{C}_{\bar u} \setminus \{0\}$ holds if and only if there exist $\varepsilon > 0$, $\tau>0$, and $\kappa >0$ such that
\begin{equation}
J''(u)v^2\geq \kappa \Vert v\Vert_{L^2(\Omega)}^2\quad \forall v\in {\mathcal C}^\tau_{\bar u} \ \text{ and }\ \forall u \in B_\varepsilon(\bar u),
\label{WWE2.18}
\end{equation}
where $B_\varepsilon(\bar u)$ denotes the ball of $L^2(\Omega)$ centered at 0 and radius $\varepsilon$.
\label{WWT2.16}
\end{theorem}

Since $C_{\bar u} \subset C^\tau_{\bar u}$ it is enough to prove that \eqref{WWE2.16} implies \eqref{WWE2.18}. The proof of this implication follows from \cite[Theorem 2.7]{Casas-Troltzsch12}. The inequality \eqref{WWE2.18} was established in that theorem for a critical cone
\[
E_{\bar u}^{\tau'} = \{v \in L^2(\Omega) : |J'(\bar u)v| \le \tau'\Vert v\Vert_{L^2(\Omega)} \text{ and  \eqref{WWE2.15}  holds}\}.
\]
The inequality on this cone implies the inequality on $C_{\bar u}^\tau$ for $\tau = \frac{\tau'}{\sqrt{|\Omega|}}$. Indeed, it is enough to observe that $C_{\bar u}^\tau \subset E_{\bar u}^{\tau'}$: for every $v \in C_{\bar u}^\tau$ we have
\[
|J'(\bar u)v|\le \int_\Omega|\nu\bar u - \bar y\bar\varphi| |v| \dx \le \tau\int_\Omega|v| \dx \le \tau\sqrt{|\Omega|}\Vert v\Vert_{L^2(\Omega)} = \tau'\Vert v\Vert_{L^2(\Omega)}.
\]

\section{Discretization and convergence} \label{S3}

In this section, we consider the numerical discretization of the control problem and prove the convergence of the discretizations. To this end we make the following additional hypothesis:

\begin{assumption} \label{A3.1}
If $\beta = +\infty$ then there exists $L_0 \in L^1(\Omega)$ such that $L(x,y) \ge L_0(x)$ for all $y \in \mathbb{R}$ and almost all $x \in \Omega$.
\end{assumption}

We discretize the problem using the finite element method. Let us consider, cf. \cite[Definition (4.4.13)]{Brenner-Scott2008}, a quasi-uniform family of triangulations $\{\mathcal{T}_{h}\}_{h>0}$ of $\bar\Omega$.
For a given triangular mesh $\mathcal{T}_h$, we define
\[Y_h = \{y_h\in C(\bar\Omega):\ y_{h|T}\in P_1(T)\ \forall T\in\mathcal{T}_h\},\]
where $P_1(T)$ is the space formed by the polynomials of degree 1 on $T$. We define the bilinear mapping $\mathfrak{a}:H^1(\Omega) \times H^1(\Omega) \longrightarrow \mathbb{R}$ given by
\[
\mathfrak{a}(y,\zeta) = \sum_{i,j=1}^{\dimension} \int_\Omega a_{i,j}(x)\partial_{x_i} y(x)\partial_{x_j}\zeta(x)\dx\quad \forall y,\zeta\in H^1(\Omega).
\]
For $u\in \A_0$, $y_h(u)$ denotes the unique solution in $Y_h$ of the equation
\begin{equation}\label{WWE3.1}\mathfrak{a}(y_h,\zeta_h) + \int_\Omega  \left(a(x,y_h)+ u y_h\right) \zeta_h \dx = \int_\Gamma g \zeta_h\dx\quad \forall \zeta_h\in Y_h.\end{equation}
The proof of the existence and uniqueness of a solution of \eqref{WWE3.1} is well
known. It is enough to apply, in a convenient way, Browder's fixed
point theorem along with the monotonicity of $a$; see, e.g., \cite[Lemma 8.2.1]{Mateos2000}.

\begin{lemma}\label{WWL3.2}
For every $u\in \A_0$, there exists a constant $c>0$, which depends on $1/\Lambda_u$ and $\Vert u\Vert_{L^2(\Omega)}$ in a monotonically non decreasing way, but it is independent of the discretization parameter $h$, such that
\begin{align}
 &\Vert y_h(u)-y_u\Vert_{L^2(\Omega)} \leq  c h^2, \label{WWE3.2}\\
&\Vert y_h(u)-y_u\Vert_{L^\infty(\Omega)} \leq  c h^{2-d/2}. \label{WWE3.3}
\end{align}
If $u\in\uad$, then $c$ is independent of $\Lambda_u$. If, further $\beta<+\infty$, $c$ is independent of $u$.
\end{lemma}
\begin{proof} We will closely follow the proof of \cite[Lemma 2 and Theorem 1]{Casas-Mateos2002}.  First, we derive the estimate,
\begin{equation} \label{WWE3.4}
\Vert y_h(u)-y_u\Vert_{H^1(\Omega)} \leq  c h \Vert y_u\Vert_{H^2(\Omega)}.
\end{equation}
For this purpose we write \eqref{WWE3.1} as follows:
\begin{equation}\label{WWE3.5} \mathfrak{a}(y_h,\zeta_h) + \int_\Omega  \left( (a_0+u) y_h + b(\cdot,y_h) \right) \zeta_h \dx = - \int_{\Omega} a(\cdot,0) \zeta_h \dx + \int_\Gamma g \zeta_h \dx\ \forall \zeta_h\in Y_h, \end{equation}
where  $b(x,y):= a(x,y)-a(x,0)-a_0(x)y.$
The above formulation is the discrete analogue of the alternative weak formulation associated to \eqref{E1.1} which is given by
\begin{equation} \label{WWE3.6}
 \mathfrak{a} (y,\zeta)  + \int_{\Omega} \left ( (a_0+u)y + b(\cdot,y) \right ) \zeta \dx = - \int_\Omega a(\cdot,0) \zeta \dx + \int_\Gamma g \zeta \dx
\end{equation}
We note that $u \in \A_0,$ $b(x,0)=0$ and $\frac{\partial b}{\partial y} (x,y) \geq 0$ due to Assumption \ref{A2.3}. Therefore, working identically to \cite[Lemma 2]{Casas-Mateos2002}, and using the fact that $y_u \in H^2(\Omega)$, we deduce \eqref{WWE3.4}.
Using \eqref{WWE3.4} and the techniques of Lemma 4 and Theorem 1 in \cite{Casas-Mateos2002} we obtain the estimates \eqref{WWE3.2} and \eqref{WWE3.3}. The dependence on $\Lambda_u$ and $\Vert u\Vert_{L^2(\Omega)}$ can be traced in the proof of \cite[Lemma 2]{Casas-Mateos2002}.
The independence of $\Lambda_u$ for $u\in\uad$ follows from Remark \ref{WWR2.10}. Finally, if $\beta<+\infty$ then $\Vert u\Vert_{L^2(\Omega)}\leq |\Omega|^{1/2}\beta$ holds and $c$ can be selected independently of $u$.
\end{proof}

\begin{corollary}  \label{WWC3.3}
For every $u\in\uad$ there exists $K_\infty>0$ that depends in a monotonically non-decreasing way on $\Vert u\Vert_{L^2(\Omega)}$ such that $\Vert y_h(u)\Vert_{L^\infty(\Omega)}\leq K_\infty$ for all $h > 0$.  If $\beta < +\infty$, $K_\infty$ can be taken independent of $u$.
\end{corollary}

\begin{proof}
  Applying the triangle inequality, estimates \eqref{WWE3.3} and \eqref{E2.3} and using that $0<\Lambda \leq\Lambda_u$, we have that
  \[\Vert y_h(u)\Vert_{L^\infty(\Omega)}\leq \Vert y_h(u)-y_u\Vert_{L^\infty(\Omega)} + \Vert y_u\Vert_{L^\infty(\Omega)} \leq c h^{2-\dimension/2} + M_{\infty},\]
  and the result follows by taking $K_\infty = c\,\text{diameter}(\Omega)^{2-\dimension/2} + M_{\infty}$.
\end{proof}

The space of discrete controls is defined by
\[
\Uh= \{u_h\in L^2(\Omega):\ u_{h|T}\in P_0(T)\ \forall T\in\mathcal{T}_h\},
\]
where $P_0(T)$ denotes the space of polynomials in $T$ of degree 0. We also set $\uadh = \Uh\cap \uad$.
We will denote $\Pi_h:L^2(\Omega) \longrightarrow \Uh$ the $L^2(\Omega)$ linear projection.
It is known that $\Pi_h u$ converges to $u$ in $L^2(\Omega)$ as $h$ tends to $0$ for all $u\in L^2(\Omega)$, and $\Pi_h u\in \uadh$ for all $u\in \uad$. Moreover, for every $u\in H^1(\Omega)$
\begin{equation}\label{WWE3.7}
\Vert u-\Pi_hu\Vert_{H^1(\Omega)^*} + h \Vert u-\Pi_hu\Vert_{L^2(\Omega)}\leq C h^2\Vert u\Vert_{H^1(\Omega)}.
\end{equation}


Our discrete problem is
\[
\Pbh \min_{u_h\in\uadh} J_h(u_h) = \frac{1}{2}\int_\Omega L(x,y_h(x))\dx + \frac{\tichonov}{2}\int_\Omega u^2_h(x)\dx.
\]

To prove convergence and obtain error estimates, we will apply the general results stated in \cite[Section 2]{Casas-Troltz2012COAP} for abstract problems. Notice that our functional $J$ is twice continuously differentiable in $L^2(\Omega)$, and hence our problem does not suffer from the two-norm discrepancy, so in all the results of \cite{Casas-Troltz2012COAP} we take $U_\infty = U_2 = L^2(\Omega)$. In the next series of lemmas we check that  the assumptions (A1)--(A6) of \cite{Casas-Troltz2012COAP} are satisfied and hence we deduce the desired strong convergence. First we observe that Assumptions (A1) and (A2) are straightforward consequences of Lemma \ref{WWL2.11} and Theorem \ref{WWT2.9}.
\begin{lemma}\label{WWL3.6}
  (Assumption (A3) in \cite{Casas-Troltz2012COAP} holds). $\uadh\subset\uad$ is convex and closed in $L^2(\Omega)$. For any $u\in\uad$ there exist elements $u_h\in\uadh$ such that $\lim_{h\to 0}\Vert u-u_h\Vert_{L^2(\Omega)}=0$.
\end{lemma}
\begin{proof}
The first statement is trivial. The second one follows taking $u_h  =\Pi_h u$.
\end{proof}
\begin{lemma}\label{WWL3.7}
  (Assumption (A4) in \cite{Casas-Troltz2012COAP} holds). For every $h>0$, if $\{u_k\}_{k=1}^\infty\subset \uadh$ and $u_k\rightharpoonup u$ in $L^2(\Omega)$, then $J_h(u)\leq \liminf_{k\to\infty} J_h(u_k)$.

  If $\beta=+\infty$, for every $h>0$ any sequence $\{u_k\}_{k=1}^\infty\subset \uadh$ such that $\lim_{k\to\infty}\Vert u_k\Vert_{L^2(\Omega)}=+\infty$ satisfies that $\lim_{k\to\infty}J_h(u_k)=+\infty$.

  If $\beta = +\infty$, any sequence $\{u_h\}_{h>0}$ such that $u_h\in\uadh$ and $\lim_{h\to 0}\Vert u_h\Vert_{L^2(\Omega)}=+\infty$ satisfies that $\lim_{h\to 0}J_h(u_h)=+\infty$.
\end{lemma}
\begin{proof}
The first statement is an immediate consequence of the fact $\uadh \subset \Uh$ and $\Uh$ is a finite dimensional space. Hence, $u_h \to u$ strongly in $L^\infty(\Omega)$ and we can easily obtain that $y_h(u_k) \to y_h(u)$ in $H^1(\Omega) \cap C(\bar\Omega)$. The second and third statements follow from Assumption \ref{A3.1}.
\end{proof}
\begin{lemma}\label{WWL3.8}
  (Assumption (A5) in \cite{Casas-Troltz2012COAP} holds). Consider a sequence $\{u_h\}_{h>0}$ such that $u_h\in\uadh$ and $u_h\rightharpoonup u$ in $L^2(\Omega)$, with $u\in\uad$. Then $J(u)\leq \liminf_{h\to 0} J_h(u_h)$ and $\lim_{h\to 0}|J_h(u_h)-J(u_h)|=0$ hold. If further the convergence $u_h \to u$ is strong in $L^2(\Omega)$, then $J(u) = \lim_{h\to 0} J_h(u_h)$ is fulfilled.
\end{lemma}
\begin{proof}
All these results are a straightforward consequence of the following property: if $u_h\rightharpoonup u$ in $L^2(\Omega)$, then $y_h(u_h)\to y_u$ in $C(\bar\Omega)$. Indeed, since $\uadh\subset\uad$ for all $h$, it is clear that $u\in\uad$. To prove the convergence we apply the triangle inequality.
\begin{align*}
  \Vert y_h(u_h)-y_u\Vert_{C(\bar\Omega)}\leq  & \Vert y_h(u_h)-y_{u_h}\Vert_{C(\bar\Omega)}+\Vert y_{u_h}-y_u\Vert_{C(\bar\Omega)}
\end{align*}
and the proof concludes applying estimate \eqref{WWE3.3} and Lemma \ref{WWL2.11}.
\end{proof}

\begin{lemma}\label{WWL3.9}
(Assumption (A6) in \cite{Casas-Troltz2012COAP} holds). If $\{u_k\}_{k=1}^\infty\subset \uad$, $u_k\rightharpoonup u$ in $L^2(\Omega)$, and $J(u_k)\to J(u)$, then $\lim_{k\to\infty}\Vert u_k-u\Vert_{L^2(\Omega)}=0$.
\end{lemma}
\begin{proof}
  From Lemma \ref{WWL2.11}, we deduce that the weak convergence $u_k\rightharpoonup u$ in $L^2(\Omega)$ implies the strong convergence $y_{u_k}\to y_u$ in $C(\bar\Omega)$, and therefore, we can pass to the limit to have
  \[\lim_{k \to \infty} \int_\Omega L(x,y_{u_k}(x))\dx  = \int_\Omega L(x,y_{u}(x))\dx.\]
  This, together with the convergence $J(u_k)\to J(u)$ implies that $\Vert u_k\Vert_{L^2(\Omega)}\to \Vert u\Vert_{L^2(\Omega)}$, which together with the weak convergence leads to the desired result.
\end{proof}
\begin{theorem}\label{WWT3.10}
  For every $h>0$ problem \Pbh has at least one solution $\bar u_h\in\uadh$. Moreover, if $\{\bar u_h\}_{h>0}$ is a sequence of solutions of \Pbh, then it is bounded in $L^2(\Omega)$. If a subsequence, denoted in the same way, $\{\bar u_h\}_{h>0}$ converges weakly to $\bar u$ in $L^2(\Omega)$ as $h\to 0$, then $\bar u$ is a solution of $\Pb$,
  \[\lim_{h\to 0} J_h(\bar u_h) = J(\bar u)\mbox{ and }\lim_{h\to 0}\Vert \bar u-\bar u_h\Vert_{L^2(\Omega)} =0.\]
\end{theorem}
\begin{proof}
This is a straightforward consequence of the previous lemmas  and \cite[Theorem 2.11]{Casas-Troltz2012COAP}
\end{proof}
The following result is a kind of reciprocal of Theorem \ref{WWT3.10}.
\begin{theorem}\label{WWT3.11}
    If $\bar u\in\uad$ is a strict local minimum of $\Pb$, then there exist $\rho > 0$ and $h_\rho >0$ such that for $h<h_\rho$ the problem
  \[
  \min_{u_h\in\bar B_\rho(\bar u)\cap \uadh} J_h(u_h)
   \]
  has a unique solution $\bar u_h$ and $\lim_{h\to 0}\Vert \bar u-\bar u_h\Vert_{L^2(\Omega)}=0$.
\end{theorem}
\begin{proof}
This follows from the previous lemmas and \cite[Theorem 2.12]{Casas-Troltz2012COAP}.
\end{proof}

\section{Error estimates for the discretization}\label{WWS4}

The aim of this section is to prove error estimates for $\bar u - \bar u_h$. In order to apply the general result \cite[Theorem 2.14]{Casas-Troltz2012COAP} we need to verify Assumption (A7) of \cite{Casas-Troltz2012COAP} about the differentiability properties of $J_h$. Let us start with the control-to-state mapping.
\begin{theorem}\label{WWT4.1}
For every $h>0$, there exists and open set $\A_h$ in $L^2(\Omega)$ such that $\A_0\subset \A_h \subset \A$ and for all $u\in\A_h$ equation \eqref{WWE3.1} has a unique solution $y_h\in Y_h$. Further, the
discrete control-to-state mapping $G_h:\A_h\longrightarrow Y_h$, $G_h(u) = y_h(u)$, is of class $C^2$. For every $u\in \A_h$ and every $v\in L^2(\Omega)$, $G'_h(u)v = z_{h,u,v}$ is the unique solution of
\begin{equation}\label{WWE4.1}
\mathfrak{a}(z_h,\zeta_h) +  \int_\Omega  [\frac{\partial a}{\partial y}(x,y_h(u))+ u] z_h \zeta_h\dx = -\int_\Omega y_h v \zeta_h\dx\ \forall \zeta_h\in Y_h.
\end{equation}
\end{theorem}
\begin{proof}
The proof is similar to the one of the continuous case (see e.g.  \cite[Theorem 2.5]{Casas-Chrysafinos-Mateos2023}) and leads to the existence of an open set $\A_h\subset L^2(\Omega)$ that may depend on $h$. For completeness, we prove the main steps. We denote by $Y_h^*$ the dual space of $Y_h$ and we introduce the functions $\mathbb{A}_h, \mathbb{F}_h:Y_h \longrightarrow Y_h^*$, $\mathcal{F}_h, \mathbb{B}_h:L^2(\Omega) \times Y_h \longrightarrow Y_h^*$, and $\mathbb{G}_h \in Y_h^*$ by
\begin{align*}
&\langle\mathbb{A}_hy_h,\zeta_h\rangle = \mathfrak{a}(y_h,\zeta_h),\ \langle\mathbb{F}_h(y_h),\zeta_h\rangle = \int_\Omega a(x,y_h)\zeta_h\dx,\\
& \langle\mathbb{B}_h(u,y_h),\zeta_h\rangle = \int_\Omega uy_h\zeta_h\dx,\ \langle\mathbb{G}_hy_h,\zeta_h\rangle = \int_\Gamma g\zeta_h\dx,\\
&\mathcal{F}_h(u,y_h) = \mathbb{A}_hy_h + \mathbb{F}_h(y_h) + \mathbb{B}_h(u,y_h) - \mathbb{G}_h.
\end{align*}
It is immediate to verify that $\mathcal{F}_h$ is of class $C^2$ and $\frac{\partial\mathcal{F}_h}{\partial y_h}(u,y_h):Y_h \longrightarrow \mathcal{L}(Y_h,Y_h^*)$ is given by
\[
\frac{\partial\mathcal{F}_h}{\partial y_h}(u,y_h)z_h = \mathbb{A}_hz_h + \mathbb{F}'_h(y_h)z_h + \mathbb{B}_h(u,z_h).
\]
Then, given $\bar u \in \A_0$ and $\bar y_h = y_h(\bar u)$, we have that $\mathcal{F}_h(\bar u,\bar y_h) = 0$ and that $\frac{\partial\mathcal{F}_h}{\partial y_h}(u,y_h)$ is an isomorphism if and only the variational problem
\[
\mathfrak{a}(z_h,\zeta_h) + \int_\Omega  \left(\frac{\partial a}{\partial y}(x,\bar y_h) + \bar u\right) z_h \zeta_h \dx = \langle\mathbb{L}_h,\zeta_h\rangle\ \forall \zeta_h\in Y_h
\]
has a unique solution $z_h \in Y_h$ for every $\mathbb{L}_h \in Y_h^*$. From the definition of $\A_0$ the coercivity of the bilinear mapping follows and, hence, the existence and uniqueness of a solution holds. Then, the implicit function theorem yields the existence of two balls $B_{\varepsilon^h_{\bar y_h}}(\bar y_h)\subset Y_h$ and $B_{\varepsilon^h_{\bar u}}(\bar u)\subset L^2(\Omega)$ such that  for every $u \in B_{\varepsilon^h_{\bar u}}(\bar u)$ there exists a unique function $y_h(u) \in B_{\varepsilon^h_{\bar y_h}}(\bar y_h)$ such that $\mathbb{F}_h(u,y_h(u)) = 0$. Arguing as in the proof of \cite[Theorem 2.5]{Casas-Chrysafinos-Mateos2023}) we deduce that $y_h(u)$ is the unique function of $Y_h$ satisfying $\mathbb{F}_h(u,y_h(u)) = 0$. Without loss of generality, we assume
\begin{equation}\label{WWE4.2}
\varepsilon^h_{\bar u}\leq\frac{\Lambda_{\bar u}}{C_{\Omega}^2},
\end{equation}
where $C_\Omega$ denotes the constant associated to the embedding $H^1(\Omega) \subset L^4(\Omega)$. Now, we define $\A_h = [\bigcup_{\bar u \in \A_0}B_{\varepsilon^h_{\bar u}}(\bar u)] \bigcap \A$. Thus we have that the mapping $G_h:\A_h \longrightarrow Y_h$ is well defined and, again by the implicit function theorem, it is of class $C^2$ and $G_h'(u)v$ is the solution of \eqref{WWE4.1}.
\end{proof}

\begin{lemma}\label{WWL4.2}
For every $h>0$, $J_h:\A_h\to\mathbb{R}$ is of class $C^2$ and for every $u\in\A_h$ and every $v,v_1,v_2\in L^2(\Omega)$
\begin{align}
 J'_h(u) v = &  \int_\Omega (\tichonov u- y_h(u)\varphi_h(u))v\dx,\label{WWE4.3}\\
 J''_{h}(u)(v_1,v_2) = & \int_\Omega \Big[ \frac{\partial^2 L}{\partial y^2}(x,y_h(u)) - \varphi_h(u)  \frac{\partial^2 a}{\partial y^2}(x,y_h(u)) \Big] z_{h,u,v_1} z_{h,u,v_2}\dx\nonumber \\
&- \int_\Omega \Big[ v_1z_{h,u,v_2} + v_2 z_{h,u,v_1} \Big] \varphi_u \dx + \tichonov \int_\Omega v_1 v_2 \dx\nonumber \\
& = \int_\Omega \left[v_2-(\varphi_h(u) z_{h,u,v_2}+y_h(u)\eta_{h,u,v_2})\right] v_1\dx \label{WWE4.4}\\
& = \int_\Omega \left[v_1-(\varphi_h(u) z_{h,u,v_1}+y_h(u)\eta_{h,u,v_1})\right] v_2\dx \ \forall u \in \A_h , \ \forall v_1, v_2 \in L^2(\Omega), \nonumber
\end{align}
where $z_{h,u,v_i} = G'_h(u)v_i$, $i=1,2$, the discrete adjoint state $\varphi_h(u)\in Y_h$ is the unique solution of
\begin{equation} \label{WWE4.5} \mathfrak{a}(\zeta_h,\varphi_h)  + \int_\Omega  [\frac{\partial a}{\partial y}(x,y_h(u))+ u] \varphi_h \zeta_h \dx = \int_\Omega \frac{\partial L}{\partial y}(x,y_h(u)) \zeta_h\dx\ \forall \zeta_h\in Y_h,
\end{equation}
and for every $v\in L^2(\Omega)$, $\eta_{h,u,v}\in Y_h$ is the unique solution of
\begin{align}
  \mathfrak{a}&(\zeta_h,\eta_h)  + \int_\Omega  [\frac{\partial a}{\partial y}(x,y_h(u))+ u] \eta_h \zeta_h \dx \nonumber \\
  = & \int_\Omega \left( \left[\frac{\partial^2 L}{\partial y^2}(x,y_h(u)) - \varphi_h(u) \frac{\partial^2 a}{\partial y^2}(x,y_h(u))\right] z_{h,u,v} -v \varphi_h(u) \right) \zeta_h\dx\ \forall \zeta_h\in Y_h\label{WWE4.6}
\end{align}
\end{lemma}
\begin{proof}
  This result is a consequence of the chain rule and Theorem \ref{WWT4.1}.
\end{proof}

\begin{lemma}\label{WWL4.3}
If $\bar u_h\in\uadh$ is a local solution of $\Pbh$, then
\begin{equation}\label{WWE4.7}
\int_\Omega (\tichonov {\bar u}_h- \bar y_h\bar \varphi_h)(u_h-\bar u_h)\dx\geq 0\quad \forall u_h\in\uadh,
\end{equation}
where $\bar y_h = y_h(\bar u_h)$ and $\bar\varphi_h = \varphi_h(\bar u_h)$.
\end{lemma}
\begin{proof}
Since $\bar u_h\in \uadh\subset \A_h$ is a local minimum of \Pbh and $J_h$ is $C^1$ around $\bar u_h$, we have
\[J'_{h}(\bar u_h)(u_h-\bar u_h)\geq 0\ \forall u_h\in\uadh,\]
and the result follows from Lemma \ref{WWL4.2}.
\end{proof}

\begin{lemma}\label{WWL4.4}
For every $u\in\uad$ we denote by $\varphi_h(u)$ the unique solution of \eqref{WWE4.5} and by $\varphi_u$ the unique solution of \eqref{WWE2.8}. Then, there exists a constant $C>0$ which depends in a monotonically non-decreasing way of $\Vert u\Vert_{L^2(\Omega)}$ but it is independent of $h$ such that
\begin{align}
&\label{WWE4.8}\Vert \varphi_u-\varphi_h(u)\Vert_{L^2(\Omega)}+ h \Vert \varphi_u-\varphi_h(u)\Vert_{H^1(\Omega)} \leq C h^2,\\
&\label{WWE4.9}\Vert \varphi_u-\varphi_h(u)\Vert_{L^\infty(\Omega)}\leq C h^{2-\dimension/2}.
\end{align}
If $\beta < +\infty$, $C$ can be chosen independently of $u$.
\end{lemma}
\begin{proof}
Let us consider $\varphi^h(u)\in H^2(\Omega)$ the unique solution of
\[
\left\{
\begin{array}{ll}
\displaystyle{A^*\varphi^h + \frac{\partial a}{\partial y}(\cdot,y_h(u))\varphi^h + u \varphi^h = \frac{\partial L}{\partial y}(\cdot, y_h(u))\mbox{ in }\Omega}\\
\partial_{\conormal_{A^*}}\varphi^h = 0\mbox{ on }\Gamma.
\end{array}
\right.
\]
Noticing that $\varphi_h(u)$ is the finite element approximation of $\varphi^h(u)$, we deduce from classical error estimates and Lemma \ref{WWL3.2}, Corollary \eqref{WWE3.3} and Assumption \eqref{WWA2.8} that
\begin{align*}
&\Vert \varphi^h(u)-\varphi_h(u)\Vert_{L^2(\Omega)}+ h \Vert \varphi^h(u)-\varphi_h(u)\Vert_{H^1(\Omega)} \leq c  h^2,\\
&\Vert \varphi^h(u)-\varphi_h(u)\Vert_{L^\infty(\Omega)}\leq c  h^{2-\dimension/2}.
\end{align*}
Denote $\tilde\varphi = \varphi_u-\varphi^h(u)$. This function satisfies the linear equation
\[
\left\{
\begin{array}{rl}
\displaystyle{A^*\tilde\varphi + \frac{\partial a}{\partial y}(\cdot,y_u)\tilde\varphi + u \tilde\varphi =} & \displaystyle\left(\frac{\partial a}{\partial y}(\cdot,y_h(u)) -\frac{\partial a}{\partial y}(\cdot,y_u)\right)\varphi^h(u) \\
 & \displaystyle{+  \frac{\partial L}{\partial y}(\cdot, y_u)- \frac{\partial L}{\partial y}(\cdot, y_h(u))\mbox{ in }\Omega}\\
\partial_{\conormal_{A^*}}\tilde\varphi = &0\mbox{ on }\Gamma.
\end{array}
\right.
\]
From the mean value theorem, we infer the existence of measurable functions $\theta_i:\Omega \longrightarrow [0,1]$, $i = 1, 2$, such that for $y_{\theta_i}= y_u + \theta_i(y_h(u) - y_u)$ we have
\[
\left\{
\begin{array}{rl}
\displaystyle{A^*\tilde\varphi + \frac{\partial a}{\partial y}(\cdot,y_u)\tilde\varphi + u \tilde\varphi =} & \displaystyle \frac{\partial a^2}{\partial y^2}(\cdot,y_{\theta_1}) (y_h(u)-y_u) \varphi^h(u) \\
 & \displaystyle{+  \frac{\partial^2 L}{\partial y^2}(\cdot, y_{\theta_2})(y_u - y_h(u))\mbox{ in }\Omega}\\
\partial_{\conormal_{A^*}}\tilde\varphi = &0\mbox{ on }\Gamma.
\end{array}
\right.
\]
Using that $0<\Lambda \leq \Lambda_u$, we deduce from Assumptions \ref{A2.2} and \ref{WWA2.8}, Corollary \ref{WWC3.3}, and estimate \eqref{E2.3} that there exists a constant $M>0$ such that $\Vert y_{\theta_1}\Vert_{L^\infty(\Omega)}\leq M$, $\Vert y_{\theta_2}\Vert_{L^\infty(\Omega)}\leq M$, $\Vert \varphi^h(u)\Vert_{L^\infty(\Omega)}\leq M$. Finally, using the techniques of \cite{Stampacchia65} we obtain from the equation satisfied by $\tilde\varphi$ that
\[
\Vert \tilde\varphi\Vert_{H^1(\Omega)} + \Vert \tilde\varphi\Vert_{L^\infty(\Omega)} \leq C (C_{L,M}+M C_{a,M}) \Vert y_h(u)-y_u\Vert_{L^2(\Omega)},
\]
and the result follows from estimate \eqref{WWE3.2}.
\end{proof}

\begin{corollary}  \label{WWC4.5}
For every $u\in\uad$ there exists $K_\infty>0$ that depends in a monotonically non-decreasing way on $\Vert u\Vert_{L^2(\Omega)}$ such that $\Vert \varphi_h(u)\Vert_{L^\infty(\Omega)}\leq K_\infty$ for every $h>0$. If $\beta < +\infty$, $K_\infty$ is independent of $u$.
\end{corollary}
\begin{proof}
  The proof is like that of Corollary \ref{WWC3.3}.
\end{proof}

In the rest of the section, $\bar u$ is a local solution that satisfies the second order sufficient condition \eqref{WWE2.16} and $\{\bar u_h\}_{h>0}$ is a sequence of local solutions of \Pbh converging to $\bar u$ in $L^2(\Omega)$. This sequence exists in virtue of Theorem \ref{WWT3.11}.
\begin{lemma}\label{WWL4.6}
(Assumption (A7) in \cite{Casas-Troltz2012COAP} holds). For  $h>0$ small enough, $J_h:\A_h\longrightarrow\mathbb{R}$ is of class $C^1$. For any $u\in\uad$ there exists $C>0$ that depends in a monotonically non-decreasing way of $\Vert u\Vert_{L^2(\Omega)}$ such that
\[|(J'_h(u)-J'(u))v|\leq C h^2 \Vert v\Vert_{L^2(\Omega)}\ \forall v\in L^2(\Omega).\]
If $\beta < +\infty$, then $C$ is independent of $u$.
\end{lemma}
\begin{proof}
The functionals $J_h$ are of class $C^1$ in $\A_h$; see Lemma \ref{WWL4.2}. Given $v\in L^2(\Omega)$, to prove the inequality we proceed as follows:
\begin{align*}
  |(J'_h(u)-J'(u))v|  = & \int_\Omega (y_u\varphi_u - y_h(u)\varphi_h(u)) v\dx\\
  = & \int_\Omega (y_u\varphi_u - y_h(u)\varphi_u  + y_h(u)\varphi_u - y_h(u)\varphi_h(u)) v\dx \\
  \leq & C  h^2 (\Vert \varphi_u\Vert_{L^\infty(\Omega)} + \Vert y_h(u)\Vert_{L^\infty(\Omega)}) \Vert v\Vert_{L^2(\Omega)},
\end{align*}
where we have used estimate \eqref{WWE3.2} and Lemma \ref{WWL4.4}. The inequality follows from  the uniform estimate for $\Vert y_h(u)\Vert_{L^\infty(\Omega)}$ established in Corollary \ref{WWC3.3}.
\end{proof}
As a consequence, we obtain the following error estimate.
\begin{theorem}\label{WWT4.7}
  Let $\bar u\in\uad$ be a local solution of \Pb satisfying the second order sufficient condition \eqref{WWE2.16} and let $\{\bar u_h\}_{h>0}$ be a sequence of local minimizers of  \Pbh converging to $\bar u$ in $L^2(\Omega)$. Then, there exists a constant $C>0$ independent of $h$ such that
  \[\Vert \bar u-\bar u_h\Vert_{L^2(\Omega)}\leq C h.\]
\end{theorem}
\begin{proof}
  From Lemmas \ref{WWL3.6}--\ref{WWL3.9} and \ref{WWL4.6}, we deduce with the help of \cite[Theorem 2.14]{Casas-Troltz2012COAP} that
  \[
  \Vert \bar u-\bar u_h\Vert_{L^2(\Omega)}\leq C \Big(h^4+\Vert \bar u-\Pi_h\bar u\Vert_{L^2(\Omega)}^2+J'(\bar u)(\Pi_h\bar u-\bar u)\Big)^{1/2}.
  \]
The result follows from the $H^1(\Omega)$ regularity of $\bar u$ obtained in Theorem \ref{WWT2.12}, the inequality
  \[
  J'(\bar u)(\Pi_h\bar u-\bar u) = \int_\Omega (\bar y\bar\varphi + \nu\bar u)(\Pi_h\bar u-\bar u)\dx\leq \Vert \bar y\bar\varphi + \nu\bar u\Vert_{H^1(\Omega)}\Vert \bar u-\Pi_h\bar u\Vert_{H^1(\Omega)^*},
  \]
and \eqref{WWE3.7}.
\end{proof}

\begin{corollary}\label{MMC4.8}
Under the assumptions of Theorem \ref{WWT4.7}, there exist constants $C>0$ and $C_\mu > 0$ independent of $h$, $C_\mu$ depending on $\mu \in (0,1)$, such that
\begin{align*}
&\Vert \bar y-\bar y_h\Vert_{L^\infty(\Omega)}+\Vert \bar \varphi-\bar \varphi_h\Vert_{L^\infty(\Omega)} \leq C h^{2-d/2},\\
&\Vert \bar u-\bar u_h\Vert_{L^\infty(\Omega)}\leq \left\{\begin{array}{ll}
  C_\mu h^\mu\ \ \forall \mu \in (0,1) & \text{ if } \dimension = 2,\\
  C h^{1/2}& \text{ if } \dimension = 3.
  \end{array}
  \right.
\end{align*}
\end{corollary}
\begin{proof}
  By the triangle inequality
  \[\Vert \bar y-\bar y_h\Vert_{L^\infty(\Omega)}\leq \Vert \bar y-y_{\bar u_h}\Vert_{L^\infty(\Omega)} + \Vert y_{\bar u_h}-\bar y_h\Vert_{L^\infty(\Omega)}\]
From Theorem \ref{WWT4.7} we have that $\{\bar u_h\}_{h>0}$ is bounded in $L^2(\Omega)$, and hence \eqref{WWE3.3} leads
\[
\Vert y_{\bar u_h}-\bar y_h\Vert_{L^\infty(\Omega)}\leq C h^{2-d/2}.
\]
Now, using that $\bar u_h \to \bar u$ in $L^2(\Omega)$ and the generalized mean value theorem we get
\begin{align*}
&\Vert \bar y-y_{\bar u_h}\Vert_{L^\infty(\Omega)} \le C_1\Vert \bar y-y_{\bar u_h}\Vert_{H^2(\Omega)} = C_1\Vert G(\bar u) - G(\bar u_h)\Vert_{H^2(\Omega)}\\
& \le C_1\max_{\theta \in [0,1]}\Vert G'(\bar u + \theta(\bar u_h - \bar u))\Vert_{\mathcal{L}(L^2(\Omega),H^2(\Omega))}\Vert \bar u - \bar u_h\Vert_{L^2(\Omega)} \le Ch.
\end{align*}
The estimate for $\bar y - \bar y_h$ follows from the previous estimates taking into account that $2-d/2\leq 1$.

For the adjoint state we proceed in a similar way. This time the generalized mean value theorem is applied to the mapping $\Phi:\A\to L^\infty(\Omega)$ given by $\Phi(u) = \varphi_u$, that is of class $C^1$. For the proof of this differentiability the reader is referred to \cite[Theorem 3.5]{Casas-Chrysafinos-Mateos2023}.

From \eqref{WWE4.7} we have that for every element $T$ and every $x\in T$
\[\bar u_h(x) = \proj_{[\alpha,\beta]}\left(\frac{1}{\tichonov \vert T \vert}\int_T \bar\varphi_h(x')\bar y_h(x')\dx'\right).\]
Using this, \eqref{WWE2.13}, and the embedding $H^2(\Omega) \subset C^{0,\mu}(\bar\Omega)$ for all $\mu \in (0,1)$ if $d=2$ and $\mu=1/2$ if $d=3$, we obtain
\begin{align*}
  |\bar u(x)- & \bar u_h(x)|\leq
   \frac{1}{\tichonov}\left\vert
  \bar\varphi(x) \bar y(x) - \frac{1}{\vert T \vert}\int_T \bar\varphi_h(x')\bar y_h(x')\dx'\right\vert \\
  = &
  \frac{1}{\tichonov}\left\vert \frac{1}{\vert T \vert}\int_T
  \left(\bar\varphi(x) \bar y(x) -  \bar\varphi_h(x')\bar y_h(x')\right)\dx'\right\vert \\
  \leq &
  \frac{1}{\tichonov\vert T \vert}
  \left(\int_T
  \left|\bar\varphi(x) \bar y(x) - \bar\varphi(x') \bar y(x')\right|\dx' +\int_T |\bar\varphi(x') \bar y(x') -\bar\varphi_h(x')\bar y_h(x')|\dx'\right)\\
  \leq &
  \frac{1}{\tichonov\vert T \vert}\left(\int_T
  \Vert \bar\varphi\bar y\Vert_{C^{0,\mu}(\bar\Omega)}  |x-x'|^\mu \dx' +\int_T \Vert \bar\varphi \bar y -\bar\varphi_h\bar y_h\Vert_{L^\infty(\Omega)}\dx'\right)\\
  \leq &
  \frac{1}{\tichonov}\left(
  \Vert \bar\varphi\bar y\Vert_{C^{0,\mu}(\bar\Omega)} h^\mu + \Vert \bar\varphi -\bar\varphi_h\Vert_{L^\infty(\Omega)} \Vert \bar y\Vert_{L^\infty(\Omega)} +\Vert \bar y-\bar y_h\Vert_{L^\infty(\Omega)}\Vert \bar\varphi_h\Vert_{L^\infty(\Omega)}\right)\\
  \leq &
  \left\{\begin{array}{cl}
  C_\mu h^\mu& \text{ if } \dimension = 2,\\
  C h^{1/2}& \text{ if } \dimension = 3,
  \end{array}
  \right.
  \end{align*}
where we have used the $L^\infty(\Omega)$ bound \eqref{E2.3}, Corollary \ref{WWC4.5}, and the obtained error estimates.
\end{proof}

\section{Superconvergence and postprocess}\label{WWS5}
The estimate of order $O(h)$ for the $L^2(\Omega)$-error in the control variable obtained in Theorem \ref{WWT4.7} is observed in numerical experiments. Using this estimate, the same order of convergence could be expected for the state and the adjoint state errors in $L^2(\Omega)$. Nevertheless, as usual in distributed control problems with regular data posed in convex domains, the experimental order of convergence obtained  for the state and the adjoint state is $O(h^2)$; see the numerical experiment in Section \ref{WWS6} below.

In \cite{Meyer-Rosch2004}, Meyer and R\"osch gave an explanation of this phenomenon in case of a control-constrained distributed linear-quadratic control problem governed by the Poisson equation with zero Dirichlet data and the control acting as a source term. The key of the proof is establishing the {\em superclosedness} property of the discrete optimal control $\bar u_h$ to the values of the optimal control $\bar u$ at the barycenters of the elements $T\in\mathcal{T}_h$. We are going to prove that our solution, also satisfies this property. This approach was first successfully used for a control problem governed by a semilinear equation in  \cite{Krumbiegel-Pfefferer2015}.

For any convex polygonal or polyhedral domain and any uniformly elliptic operator with Lipschitz coefficients $A$, there exists $p_{A,\Omega}>2$ such that for every $r\in (2,p_{A,\Omega})$ and every $f\in L^r(\Omega)$, the problem
        \[A y = f\text{ in }\Omega,\ \partial_n y = 0 \text{ on }\Gamma\]
        has a unique solution in $W^{2,r}(\Omega)/\mathbb{R}$. For the proof of this result the reader is referred to \cite[Theorem 4.4.3.7]{Grisvard85} for $\dimension =2$ and to \cite[Corollary 3.12]{Dauge1992} for $\dimension = 3$; see also \cite[Theorem 1.10, page 15]{Girault-Raviart86}.
Through this section we will suppose that the following assumption holds.

\begin{assumption}\label{WWA2.13}
We assume that $p_{A,\Omega}>\dimension$ and there exists some $p \in (\dimension,p_\Omega)$ such that $a(\cdot,0) \in L^p(\Omega)$, the function $y_g$ introduced in Assumption \ref{A2.3} satisfies $y_g \in W^{2,p}(\Omega)$, and the functions $L_M$ in Assumption \ref{WWA2.8} belongs to $L^p(\Omega)$.
\end{assumption}

\begin{proposition}\label{WWC2.14}
  Suppose  Assumption \ref{WWA2.13} holds. If $\bar u\in\uad$ is a local minimizer of \Pb then $\bar y,\bar\varphi\in W^{2,p}(\Omega)$ and $\bar u\in  H^1(\Omega) \cap C^{0,1} (\bar \Omega)$.
\end{proposition}
\begin{proof}
Since $\bar u\bar y$ is a continuous function, we have that $A\bar y\in L^p(\Omega)$. Using that $\partial_{\conormal_A} \bar y$ is the normal derivative of a function in $W^{2,p}(\Omega)$ and exploiting the homogeneous equation satisfied by $y_u-y_g$  as in our proof of Theorem \ref{WWT2.5}, we deduce from the  results in \cite[Theorem 4.4.3.7]{Grisvard85} and \cite{Dauge1992} that $\bar y\in W^{2,p}(\Omega)$.

The same kind of reasoning leads to $\bar\varphi\in W^{2,p}(\Omega)$. Therefore $\bar y\bar\varphi\in W^{2,p}(\Omega)\hookrightarrow C^{0,1}(\bar\Omega)$ and the result follows from \eqref{WWE2.13}.
\end{proof}

Through this section we will suppose that $\bar u \in C^{0,1}(\bar\Omega)$ is a local solution of \Pb satisfying the second order sufficient condition \eqref{WWE2.16} and $\{\bar u_h\}_{h>0}$ is a sequence of local minimizers of  \Pbh converging to $\bar u$ in $L^2(\Omega)$. We will also denote $\bar y,\bar\varphi\in W^{2,p}(\Omega)$ the state and adjoint state related to $\bar u$ and $\bar y_h,\bar\varphi_h\in Y_h$ the discrete state and adjoint state related to $\bar u_h$. As in \cite{Krumbiegel-Pfefferer2015}, for every $h>0$ we denote
\begin{align*}
  \mathcal{T}_1 = & \big\{T:\exists x_a,x_i\in T\mbox{ such that }\bar u(x_i)\in (\alpha,\beta)\mbox{ and }\bar u(x_a)\in\{\alpha,\beta\}\big\} \\
  \mathcal{T}_2 = & \mathcal{T}_h\setminus \mathcal{T}_1,
\end{align*}
and for any $\mathcal{T}\subset\mathcal{T}_h$ we define the set $\Omega_{\mathcal{T}} = \cup\{T\in \mathcal{T}\}.$
\begin{assumption}\label{WWA5.1} There exists $C>0$ independent of $h$ such that
  \[|\Omega_{\mathcal{T}_1}|\leq C h.\]
\end{assumption}
For every $T\in\mathcal{T}_h$, we set $x_T$ as the barycenter of $T$ if $T\in \mathcal{T}_2$ and $x_T=x_a$ if $T\in \mathcal{T}_1$, for any fixed $x_a\in T$ such that $\bar u(x_a)\in\{\alpha,\beta\}$.
Notice that, using the Lipschtiz property of $\bar u$, for $h>0$ small enough $\bar u(x)$ cannot achieve both values $\alpha$ and $\beta$ in the same element $T$.
Finally, we define $w_h\in\uadh$ by $w_h(x) = \bar u(x_T)$ if $x\in T$.

We first prove some approximation properties of $w_h$.
\begin{lemma}\label{WWL5.3}Suppose assumptions \ref{WWA5.1} and \ref{WWA2.13} hold. Then, there exists $C>0$ independent of $h$ such that,
\begin{align}
&\Vert \bar u- w_h\Vert_{L^\infty(\Omega)} \leq  L_{\bar u} h,     \label{WWE5.1}\\
&\Vert \Pi_h \bar u- w_h\Vert_{L^1(\Omega)} \leq C h^2, \label{WWE5.2}
\end{align}
where $L_{\bar u}$ is the Lipschitz constant of $\bar u$. Moreover, there exist $\tau_0 > 0$ and $h_0 > 0$ such that $\bar u_h-w_h\in \mathcal{C}^\tau_{\bar u}$ for all $h \le h_0$ and $\tau \in (0,\tau_0)$.
\end{lemma}
\begin{proof}
 Let us take $T \in \mathcal{T}_h$ arbitrary and $x \in T$, then we have
  \[
  |\bar u(x) - w_h(x)| = |\bar u(x) - \bar u(x_T)| \le L_{\bar u}|x - x_T| \le L_{\bar u}h,
  \]
   which proves \eqref{WWE5.1}. Let us prove \eqref{WWE5.2}.
  First, noting that $w_h$ is constant in every element $T$
  \begin{align*}
  \int_\Omega |\Pi_h\bar u-w_h|\dx & = \sum_{T\in\mathcal{T}_h}\int_T |\Pi_h (\bar u-w_h)| \dx =
  \sum_{T\in\mathcal{T}_h}\int_T \left|\frac{1}{|T|} \int_T (\bar u-w_h)\mathrm{d}x' \right| \dx\\
  = & \sum_{T\in\mathcal{T}_h} \left| \int_T (\bar u-w_h)\dx' \right|.
  \end{align*}
  For every $T\in \mathcal{T}_2$, either $\bar u = \frac{1}{\tichonov}\bar y\bar\varphi\in W^{2,p}(T)\hookrightarrow H^2(T)$ or $\bar u$ is constant. Using \cite[Lemma 3.2]{Meyer-Rosch2004} and the discrete Cauchy-Scwhartz inequality, we obtain
  \begin{align*}
    \sum_{T\in \mathcal{T}_2} \left| \int_T (\bar u-w_h)\dx \right| \leq &  \sum_{T\in \mathcal{T}_2} c h^2 \sqrt{|T|} |\bar u|_{H^2(T)} \leq c h^2 \sqrt{|\Omega_{\mathcal{T}_2}|} |\bar u|_{H^2(\Omega_{\mathcal{T}_2})} \le \frac{ch^2}{\nu}\sqrt{|\Omega|}|\bar y\bar\varphi|_{H^2(\Omega)}.
  \end{align*}
  For the other part, we use the Lispchitz continuity of $\bar u$ and Assumption \ref{WWA5.1}.
  \begin{align*}
    \sum_{T\in \mathcal{T}_1} \left| \int_T (\bar u-w_h)\dx \right| \leq &  \sum_{T\in \mathcal{T}_1}  \int_T |\bar u-\bar u(x_T)|\dx \leq L_{\bar u}\sum_{T\in \mathcal{T}_1}  \int_T |x-x_T|\dx\\
\leq&    L_{\bar u}\sum_{T\in \mathcal{T}_1} h |T| \leq L_{\bar u}h |\Omega_{\mathcal{T}_1}| \leq CL_{\bar u} h^2.
  \end{align*}
This concludes the proof of \eqref{WWE5.2}.

For the proof of  $\bar u_h-w_h\in \mathcal{C}^\tau_{\bar u}$ for $h$ and $\tau$ small enough, the reader is referred to the proof of \cite[Lemma 6]{Krumbiegel-Pfefferer2015}. Notice that we can follow it thanks to the uniform convergence proved in Corollary \ref{MMC4.8}.
\end{proof}

Next, we establish the approximation properties of the discrete states and adjoint states related to $w_h$.
\begin{lemma}\label{WWL5.4}
Suppose assumptions \ref{WWA5.1} and \ref{WWA2.13} hold. Then, there exists $C>0$ independent of $h$ such that,
\begin{align}
&\Vert \bar y- y_h(w_h)\Vert_{L^2(\Omega)} \leq C h^2, \label{WWE5.3}\\
&\Vert \bar\varphi- \varphi_h(w_h)\Vert_{L^2(\Omega)} \leq C h^2. \label{WWE5.4}
\end{align}
\end{lemma}
\begin{proof}
  Taking into account estimate \eqref{WWE3.2}, estimate \eqref{WWE5.3} follows from
  \begin{equation}\label{WWE5.5}
    \Vert y_h(\bar u)- y_h(w_h)\Vert_{L^2(\Omega)}  \leq C h^2.
  \end{equation}
Let us check \eqref{WWE5.5}. Define
\[
\delta_h(x) = \left\{\begin{array}{cl}
                         \displaystyle\frac{a(x,y_h(\bar u)(x)) - a(x,y_h(w_h)(x))}{y_h(\bar u)(x) - y_h(w_h)(x)} + \bar u(x) & \text{if } y_h(\bar u)(x) - y_h(w_h)(x)\neq 0,\vspace{2mm}\\
                        a_0(x)+\alpha & \text{otherwise}.
                       \end{array}\right.
\]
From Corollary \ref{WWC3.3}, \eqref{WWE5.1}, Assumption \ref{A2.2}, and Assumption \ref{WWA2.7}, we deduce that $\delta_h\in L^\infty(\Omega)$ and there exists a constant $D_\infty$ such that $ 0\le a_0(x)+\alpha \leq \delta_h(x)\leq D_\infty$ for all $x\in \Omega$ and all $h>0$.
  Hence, there exist unique $\phi\in H^2(\Omega)$ and $\phi_h\in Y_h$ such that
\begin{align}
&\mathfrak{a}(\zeta,\phi) + \int_\Omega \delta_h  \phi \zeta \dx  = \int_\Omega (y_h(\bar u)-y_h(w_h)) \zeta \dx\quad\forall \zeta \in H^1(\Omega), \label{WWE5.6} \\
&\mathfrak{a}(\zeta_h,\phi_h) + \int_\Omega \delta_h  \phi_h \zeta_h \dx = \int_\Omega (y_h(\bar u)-y_h(w_h)) \zeta_h\dx\quad \forall \zeta_h\in Y_h. \label{WWE5.7}
  \end{align}
  Notice that $\phi_h$ is the finite element approximation of $\phi$. We have that
  \begin{equation}\label{WWE5.8}
  \Vert \phi\Vert_{H^2(\Omega)}\leq C \Vert y_h(\bar u)-y_h(w_h)\Vert_{L^2(\Omega)},
  \end{equation}
  \begin{equation}  \label{WWE5.9}
  \Vert \phi-\phi_h\Vert_{L^2(\Omega)} + h \Vert \phi-\phi_h\Vert_{H^1(\Omega)} + h^{\dimension/2}\Vert \phi-\phi_h\Vert_{L^\infty(\Omega)}  \leq C h^2 \Vert y_h(\bar u)-y_h(w_h)\Vert_{L^2(\Omega)}.
  \end{equation}
As in Corollary \ref{WWC3.3}, estimates \eqref{WWE5.8} and \eqref{WWE5.9} and the embedding $H^2(\Omega)\hookrightarrow L^\infty(\Omega)$ imply the existence of $C>0$  independent of $h$ such that
  \begin{equation}\label{WWE5.10}\Vert \phi_h\Vert_{L^\infty(\Omega)}+\Vert \phi_h\Vert_{H^1(\Omega)}\leq C\Vert y_h(\bar u)-y_h(w_h)\Vert_{L^2(\Omega)}\end{equation}
  for all $h>0$.
Testing equation \eqref{WWE5.7} with $\zeta_h = y_h(\bar u)-y_h(w_h)$ and using the definition of $\delta_h $ and the nonlinear finite dimensional equations satisfied by $y_h(\bar u)$ and $y_h(w_h)$, we obtain
  \begin{align}
    \Vert  y_h(\bar u)& -y_h(w_h)\Vert_{L^2(\Omega)}^2 =   \mathfrak{a}(y_h(\bar u)-y_h(w_h),\phi_h)\notag
    \\
    & + \int_\Omega (a(\cdot,y_h(\bar u)) - a(\cdot,y_h(w_h)))\phi_h\dx\notag
     +\int_\Omega \bar u(x)(y_h(\bar u) - y_h(w_h))\phi_h\dx\notag
    \\
    = & \mathfrak{a}(y_h(\bar u),\phi_h) + \int_\Omega (a(\cdot,y_h(\bar u))+\bar u y_h( \bar u)) \phi_h\dx\notag
    \\
     & -\Big[\mathfrak{a}(y_h(w_h),\phi_h) + \int_\Omega (a(\cdot,y_h(w_h))+ w_h(x)y_h(w_h)) \phi_h\dx\Big] \notag
     \\
     & + \int_\Omega(w_h(x)-\bar u(x))y_h(w_h)(x)\phi_h(x)\dx\notag
     \\
  = &\int_\Omega(w_h(x)-\bar u(x))y_h(w_h)(x)\phi_h(x)\dx.\label{WWE5.11}
  \end{align}
Using Corollary \ref{WWC3.3}, \eqref{WWE3.7}, \eqref{WWE5.2}, and \eqref{WWE5.10} we estimate
\begin{align}
\int_\Omega & (w_h(x)-\bar u(x))  y_h(w_h)(x)\phi_h(x)\dx \notag \\
=  &\int_\Omega  (w_h(x)-\Pi_h\bar u(x))y_h(w_h)(x)\phi_h(x)\dx\notag + \int_\Omega(\Pi_h\bar u(x)-\bar u(x))y_h(w_h)(x)\phi_h(x)\dx  \notag \\
 \leq &\  K_\infty \Vert w_h-\Pi_h \bar u\Vert_{L^1(\Omega)} \Vert \phi_h\Vert_{L^\infty(\Omega)}  + \Vert \Pi_h\bar u-\bar u\Vert_{(H^1(\Omega)^*)}  \Vert y_h(w_h) \phi_h\Vert_{H^1(\Omega)} \notag\\
\leq &\ C h^2   \Vert y_h(\bar u)-y_h(w_h)\Vert_{L^2(\Omega)}. \label{WWE5.12}
\end{align}
Now \eqref{WWE5.5} follows from \eqref{WWE5.11} and \eqref{WWE5.12}.

The proof of \eqref{WWE5.4} follows the same lines. Since the adjoint state equation is linear, the appropriate reaction coefficient is
\[
\delta_h  = \frac{\partial a}{\partial y}(x,y_h(\bar u)) + \bar u,
\]
and the auxiliary adjoint equation is
\[
\mathfrak{a}(\eta,\zeta) + \int_\Omega \delta_h  \eta \zeta \dx  = \int_\Omega (\varphi_h(\bar u)-\varphi_h(w_h)) \zeta \dx\quad \forall \zeta \in H^1(\Omega).
\]
Appropriate estimates of the terms
 \begin{align}
   &\int_\Omega \left(\frac{\partial a}{\partial y}(x,y_h(\bar u))-\frac{\partial a}{\partial y}(x,y_h(w_h))\right)\varphi_h(w_h)\eta_h\dx,
   \label{WWE5.13}\\
   &\int_\Omega \left(\frac{\partial L}{\partial y}(x,y_h(\bar u))-\frac{\partial L}{\partial y}(x,y_h(w_h))\right)\eta_h\dx,
   \label{WWE5.14}
\end{align}
where $\eta_h$ is the finite element approximation of $\eta$, that are needed in the process, are obtained using \eqref{WWE5.3} and Assumptions \ref{A2.2} and \ref{WWA2.8}.
\end{proof}

Next we prove a discrete variational inequality satisfied by $w_h$, which will be used to compare it with $\bar u_h$ later.
\begin{lemma}\label{WWL5.5} Let us take $T \in \mathcal{T}_h$ arbitrary.
 Then we have $(\nu w_h(x_T) - \bar y(x_T)\bar\varphi(x_T))(\bar u_h(x_T) - w_h(x_T)) \ge 0$.
\end{lemma}
\begin{proof}
From \eqref{WWE2.13} we deduce that $(\nu \bar u(x) - \bar y(x)\bar\varphi(x))(\xi - \bar u(x)) \ge 0$ $\forall \xi \in [\alpha,\beta]$. Since $\bar u_h(x_T)\in[\alpha,\beta]$, this implies that $(\nu\bar u(x_T)- \bar y(x_T)\bar\varphi(x_T))(\bar u_h(x_T)-\bar u(x_T))\geq 0$, which can be written as $(\nu w_h(x_T)-\bar y(x_T)\bar\varphi(x_T))(\bar u_h(x_T)- w_h(x_T))\geq 0$.
\end{proof}

\begin{lemma}\label{ZZL5.8}
  Suppose Assumption \ref{WWA5.1} holds. Then, there exists $C$ independent of $h$ such that the following inequalities hold
  \[\Vert \bar u_h-w_h\Vert_{L^2(\Omega_{\mathcal{T}_1})}\leq C h^{3/2}\text{ and }
  \Vert \bar u_h-w_h\Vert_{L^2(\Omega_{\mathcal{T}_2})}\leq C h^{2}.\]
\end{lemma}
\begin{proof}
Given $i \in \{1,2\}$, we define $\tilde w_h\in\uadh$ by $\tilde w_{h\mid_T} = w_h$ if $T\in \mathcal{T}_i$ and $\tilde w_{h\mid_T} = \bar u_h$ if $T\in \mathcal{T}_{h} \setminus \mathcal{T}_i$. Taking $u_h = \tilde w_h$ in \eqref{WWE4.7} we get
\[
\int_{\Omega_{\mathcal{T}_i}}((\nu\bar u_h-\bar y_h\bar\varphi_h)(w_h-\bar u_h)\dx\geq 0.
\]
Let $d_h\in \Uh$ be defined as $d_h(x) = \bar y(x_T)\bar\varphi(x_T)$ if $x\in T$. The above inequality and Lemma \ref{WWL5.5} imply
\[
\int_{\Omega_{\mathcal{T}_i}}(\bar y_h\bar\varphi_h-d_h+\nu(w_h-\bar u_h))(\bar u_h - w_h)\dx\geq 0.
\]
This leads to
\begin{align*}
 0 &\leq \int_{\Omega_{\mathcal{T}_i}}(\bar y\bar\varphi-d_h)(\bar u_h-w_h)\dx
+ \int_{\Omega_{\mathcal{T}_i}}[\bar y_h(w_h)\bar\varphi_h(w_h)-\bar y\bar\varphi](\bar u_h-w_h)\dx \\
& +\int_{\Omega_{\mathcal{T}_i}}\{[\nu w_h-\bar y_h(w_h)\bar\varphi_h(w_h)] - [\nu\bar u_h-\bar y_h\bar\varphi_h]\}(\bar u_h-w_h)\dx = I + II + III
\end{align*}
Let us estimate $I$ in the two cases $i=1$ or $i=2$:

 If $T\in \mathcal{T}_2$, then $x_T$ is the barycenter of $T$. Hence, using the fact that $\bar y\bar\varphi\in H^2(T)$ and that $\bar u_h-w_h$ is constant on $T$, we infer from Lemma \ref{WWL5.5} that
\[
\int_T(\bar y\bar\varphi(\bar u_h-w_h)-d_h(\bar u_h-w_h))\dx\leq c h^2 \left( |\bar u_h-w_h| \sqrt{|T|}\right) |\bar y\bar\varphi|_{H^2(T)}.
\]
  Adding up over all $T\in \mathcal{T}_2$ and applying the discrete Cauchy-Schwartz inequality, we have
  \begin{align}
  \sum_{T\in \mathcal{T}_2} \int_T(\bar y\bar\varphi(\bar u_h-w_h)-d_h(\bar u_h-w_h))\dx \leq & c h^2 \Vert \bar u_h-w_h\Vert_{L^2(\Omega_{\mathcal{T}_2})} |\bar y\bar\varphi|_{H^2(\Omega_{\mathcal{T}_2})}\label{WWE5.16}\\
  \leq & c h^2 \Vert \bar u_h-w_h\Vert_{L^2(\Omega)} |\bar y\bar\varphi|_{H^2(\Omega)}.\notag
  \end{align}
  If  $T\in \mathcal{T}_1$, then $x_T$ is not the barycenter of $T$ and the numerical integration formula is not of order $h^2$, but $h$. From Assumption \ref{WWA2.13}, we have that $\bar y\bar\varphi\in W^{2,p}(\Omega)\hookrightarrow C^{0,1}(\bar\Omega)$. Using this and Assumption \ref{WWA5.1}.
  \[\int_T(\bar y\bar\varphi(\bar u_h-w_h)-d_h(\bar u_h-w_h))\dx \leq C h|\bar u_h-w_h| |T| \Vert \bar y\bar\varphi\Vert_{C^{0,1}(\bar\Omega)}\]
Adding up over all $T\in \mathcal{T}_1$ and applying the discrete Cauchy-Schwartz inequality and Assumption \ref{WWA5.1}, we obtain
  \begin{align*}
  \sum_{T\in \mathcal{T}_1} \int_T(\bar y\bar\varphi(\bar u_h-w_h)-d_h(\bar u_h-w_h))\dx \leq & c h \Vert \bar u_h-w_h\Vert_{L^2(\Omega_{\mathcal{T}_1})} \sqrt{|\Omega_{\mathcal{T}_1}|}\\
  \leq & c h^{3/2} \Vert \bar u_h-w_h\Vert_{L^2(\Omega_{\mathcal{T}_1})},
  \end{align*}

 From Lemma \ref{WWL5.4}, we deduce that that $II\leq C h^2\Vert \bar u_h-w_h\Vert_{L^2(\Omega_{\mathcal{T}_i})}$ for $i=1, 2$.

Using the mean value theorem, the term III can be written as
\begin{align}
III = & (J'_h(w_h)-J'_h(\bar u_h))\chi_{\Omega_{\mathcal{T}_i}}(\bar u_h-w_h) = - J''_h(u_{\theta_h})( \bar u_h-w_h, \chi_{\Omega_{\mathcal{T}_i}}(\bar u_h-w_h))\notag \\
= & -J''_h(u_{\theta_h})( \chi_{\Omega_{\mathcal{T}_i}}(\bar u_h-w_h), \chi_{\Omega_{\mathcal{T}_i}}(\bar u_h-w_h)) \notag \\
  & -J''_h(u_{\theta_h})( \chi_{\Omega_{\mathcal{T}_j}}(\bar u_h-w_h), \chi_{\Omega_{\mathcal{T}_i}}(\bar u_h-w_h))\label{ZZE5.17},
\end{align}
where $u_{\theta_h}(x) = \bar u_h(x) + \theta_h(w_h(x) - \bar u_h(x))$ and $j \in \{1, 2\} \setminus \{i\}$. The first term in the right hand side of \eqref{ZZE5.17} can be estimated as follows.
Since $\bar u_h-w_h\in {\mathcal C}^\tau_{\bar u}$, so does $\chi_{\Omega_{\mathcal{T}_j}}(\bar u_h-w_h)$. Using the second order sufficient optimality condition \eqref{WWE2.18}, Theorems \ref{WWT2.16} and  \ref{WWTA.9}, writing $J''_h(u_\theta) = (J''_h(u_\theta)-J''(u_\theta)) + J''(u_\theta)$, we deduce that for $h>0$ small enough
\begin{equation}\label{ZZE5.15}
-J''_h(u_{\theta_h})( \chi_{\Omega_{\mathcal{T}_i}}(\bar u_h-w_h), \chi_{\Omega_{\mathcal{T}_i}} (\bar u_h-w_h)) \leq -\frac{\kappa}{2}\Vert \bar u_h-w_h\Vert_{L^2(\Omega_{\mathcal{T}_i})}^2.
\end{equation}
For the second term, denoting $v_k = \chi_{\Omega_{\mathcal{T}_k}}(\bar u_h-w_h)$ for $k=1, 2$ and using \eqref{WWE4.4}, the fact that $\mathcal{T}_1\cap \mathcal{T}_2 = \emptyset$, and Corollaries \ref{WWC3.3}, \ref{WWC4.5}, \ref{WWCA.5}, and the estimate $\Vert \eta_{h,u,v}\Vert_{H^1(\Omega)} \le c\Vert v\Vert_{L^2(\Omega)}$ for all $u \in \uad$ and $v \in L^2(\Omega)$, we have
\begin{align*}
  -J''_h(u_{\theta_h})( \chi_{\Omega_{\mathcal{T}_j}}(\bar u_h-w_h), & \chi_{\Omega_{\mathcal{T}_i}}(\bar u_h-w_h)) \\ = &
  -\int_\Omega (\varphi_h(u_{\theta_h})z_{h,u_{\theta_h},v_2}+y_h(u_{\theta_h})\eta_{h,u_{\theta_h},v_2})v_1\dx  \\
  &-\int_\Omega \chi_{\Omega_{\mathcal{T}_1}}(\bar u_h-w_h) \chi_{\Omega_{\mathcal{T}_2}}(\bar u_h-w_h)\dx\\
  \leq\ & 2K_\infty K_1c \Vert \bar u_h-w_h\Vert_{L^2(\Omega_{\mathcal{T}_2})} \Vert \bar u_h-w_h\Vert_{L^1(\Omega_{\mathcal{T}_1})}
\end{align*}
Applying the continuous and discrete versions of Cauchy-Schwarz inequality and Assumption \ref{WWA5.1}, we have
  \begin{align*}
 \Vert \bar u_h-w_h\Vert_{L^1(\Omega_{\mathcal{T}_1})} = & \sum_{T\in \mathcal{T}_1}\int|\bar u_h-w_h|\dx\leq \sum_{T\in \mathcal{T}_1} |T|^{1/2} \Vert \bar u_h-w_h\Vert_{L^2(T)}\\
 \leq & |\Omega_{\mathcal{T}_1}|^{1/2} \Vert \bar u_h-w_h\Vert_{L^2(\Omega_{\mathcal{T}_1})}\leq C h^{1/2}\Vert \bar u_h-w_h\Vert_{L^2(\Omega_{\mathcal{T}_1})}
 \end{align*}
From the last two estimates we gather the existence of $C>0$ such that
\begin{equation}\label{ZZE5.18}
  -J''_h(u_{\theta_h})( \chi_{\Omega_{\mathcal{T}_j}}(\bar u_h-w_h), \chi_{\Omega_{\mathcal{T}_i}}(\bar u_h-w_h)) \leq
   C h^{1/2} \Vert \bar u_h-w_h\Vert_{L^2(\Omega_{\mathcal{T}_2})}  \Vert \bar u_h-w_h\Vert_{L^2(\mathcal{T}_1)}
\end{equation}
From \eqref{ZZE5.17}--\eqref{ZZE5.18} we obtain:
\begin{equation}\label{ZZE5.19}III\leq -\frac{\kappa}{2}\Vert \bar u_h-w_h\Vert_{L^2(\Omega_{\mathcal{T}_2})}^2+ Ch^{1/2}\Vert \bar u_h-w_h\Vert_{L^2(\Omega_{\mathcal{T}_1})} \Vert \bar u_h-w_h\Vert_{L^2(\Omega_{\mathcal{T}_2})}.\end{equation}

Gathering the estimates for the three terms I, II and III in the case $i=1$ we have that
\[
\frac{\kappa}{2}\Vert \bar u_h-w_h\Vert_{L^2(\Omega_{\mathcal{T}_1})}^2 \leq Ch^{1/2}\Vert \bar u_h-w_h\Vert_{L^2(\Omega_{\mathcal{T}_1})} \Vert \bar u_h-w_h\Vert_{L^2(\Omega_{\mathcal{T}_2})} + c h^{3/2} \Vert \bar u_h-w_h\Vert_{L^2(\Omega_{\mathcal{T}_1})}
\]
which leads to
\begin{equation}\label{ZZE5.20}
  \frac{\kappa}{2}\Vert \bar u_h-w_h\Vert_{L^2(\Omega_{\mathcal{T}_1})} \leq Ch^{1/2}\Vert \bar u_h-w_h\Vert_{L^2(\Omega_{\mathcal{T}_2})} + c h^{3/2}.
\end{equation}
In a similar way, in the case $i=2$ we obtain
\begin{equation}\label{ZZE5.21}
  \frac{\kappa}{2}\Vert \bar u_h-w_h\Vert_{L^2(\Omega_{\mathcal{T}_2})} \leq Ch^{1/2}\Vert \bar u_h-w_h\Vert_{L^2(\Omega_{\mathcal{T}_1})} + c h^2.
\end{equation}
Inserting \eqref{ZZE5.20} in \eqref{ZZE5.21} we obtain
\[  \frac{\kappa}{2}\Vert \bar u_h-w_h\Vert_{L^2(\Omega_{\mathcal{T}_2})} \leq C h\Vert \bar u_h-w_h\Vert_{L^2(\Omega_{\mathcal{T}_2})} + ch^2  \]
and the estimate for $\Vert \bar u_h-w_h\Vert_{L^2(\Omega_{\mathcal{T}_2})}$ follows for all $h$ such that $Ch < \kappa/4$.
The estimate for $\Vert \bar u_h-w_h\Vert_{L^2(\Omega_{\mathcal{T}_1})}$ now follows from this one and \eqref{ZZE5.20}.
\end{proof}

\begin{theorem}\label{WWT5.9}
  Suppose assumptions \ref{WWA5.1} and \ref{WWA2.13} hold. Then there exists $C>0$ such that
\begin{align}
&\Vert \bar y_h-\bar y\Vert_{L^2(\Omega)}\leq C h^2,\label{WWE5.19}\\
&\Vert \bar \varphi_h-\bar \varphi\Vert_{L^2(\Omega)}\leq C h^2\label{WWE5.20}.
&\end{align}
\end{theorem}
\begin{proof}
We prove that $\Vert \bar y_h-y_h(w_h)\Vert_{L^2(\Omega)} = O(h^2)$, then \eqref{WWE5.19}  follows from Lemma \ref{WWL5.4}.  First we demonstrate the inequality
\begin{equation}\label{WWE5.21}
\Vert \bar y_h-y_h(w_h)\Vert_{L^2(\Omega)}\leq C \Vert \bar u_h-w_h\Vert_{L^1(\Omega)}.
\end{equation}
   We repeat the first steps of the proof  of Lemma \ref{WWL5.4} replacing $\bar u$ by $\bar u_h$ in the definitions of $\delta_h$, $\phi$ and $\phi_h$ in \eqref{WWE5.6} - \eqref{WWE5.7} to obtain
  \[ \Vert \bar y_h-y_h(w_h)\Vert_{L^2(\Omega)}^2\leq \int_\Omega (w_h-\bar u_h) y_h(w_h) \phi_h\dx \leq \Vert y_h(w_h)\Vert_{L^\infty(\Omega)} \Vert \phi_h\Vert_{L^\infty(\Omega)} \Vert \bar u_h-w_h\Vert_{L^1(\Omega)}.\]
  Finally, using standard uniform error estimates for linear equations we obtain
  \[\Vert \phi_h\Vert_{L^\infty(\Omega)} \leq \Vert \phi_h-\phi\Vert_{L^\infty(\Omega)} + \Vert \phi\Vert_{L^\infty(\Omega)}\leq (Ch^{2-\dimension/2}+1)\Vert \bar y_h-y_h(w_h)\Vert_{L^2(\Omega)}.\]
  This estimate, together with Corollary \ref{WWC3.3}, leads to estimate \eqref{WWE5.21}.

  Now, applying Cauchy-Schwarz inequality, Assumption \ref{WWA5.1}, Lemma \ref{ZZL5.8} we infer
  \begin{align*}
    \Vert \bar u_h-w_h\Vert_{L^1(\Omega)} \leq & \Vert \bar u_h-w_h\Vert_{L^1(\Omega_{\mathcal{T}_1})} + \Vert \bar u_h-w_h\Vert_{L^1(\Omega_{\mathcal{T}_2})} \\
    \leq  & |\Omega_{\mathcal{T}_1}|^{1/2} \Vert \bar u_h-w_h\Vert_{L^2(\Omega_{\mathcal{T}_1})} +
    |\Omega_{\mathcal{T}_2}|^{1/2}
    \Vert \bar u_h-w_h\Vert_{L^2(\Omega_{\mathcal{T}_2})}  \leq  C h^2.
  \end{align*}
The proof of \eqref{WWE5.20} follows the same lines. Terms analogous to those of \eqref{WWE5.13} and \eqref{WWE5.14} with $\bar u$ replaced by $\bar u_h$ appear, and can be bounded with the help of Assumptions \ref{A2.2} and \ref{WWA2.8} and the estimate  $\Vert \bar y_h-y_h(w_h)\Vert_{L^2(\Omega)}\leq C h^2$ just obtained.
\end{proof}

\begin{remark}
The postprocessing approach suggested by Meyer and R\"osch \cite{Meyer-Rosch2004} consists of setting
\[
\tilde u_h(x) = \proj_{[\alpha,\beta]}\left(\frac{1}{\nu}\bar y_h(x)\bar\varphi_h(x)\right).
\]
Then, as a consequence of Theorem \ref{WWT5.9} we obtain he error estimate:
\[
\Vert \tilde u_h(x)-\bar u\Vert_{L^2(\Omega)} \le Ch^2.
\]
\end{remark}

\section{A numerical experiment}\label{WWS6}
We solve the following problem, introduced in \cite{Casas-Chrysafinos-Mateos2023}.
Let $\Omega$ be $(0,1)^2\subset\mathbb{R}^2$, $A=-\Delta$, $g(x)=0$,
\[f(x,y) = y^3|y|+2y-100\sin(2\pi x_1)\sin(\pi x_2),\]
$\tichonov = 0.05$, $\umin = -1$, $\umax = 1$, $L(x,y) = 0.5(y-y_d(x))$, where
\[y_d(x) = -64 x_1(1-x_1) x_2(1-x_2).\]
The discretization is done using a family of quasi-uniform meshes obtained by diadic refinement, with sizes $h_j = 2^{-j}\sqrt{2}$, $j\in\{3, \ldots,11\}$, the finest mesh having $4.2\times 10^6$ nodes and $8.4\times 10^6$ elements. Each problem is solved using the semismooth Newton method described in \cite{Casas-Chrysafinos-Mateos2023}, where its convergence analysis is also studied. We measure the error in $L^2(\Omega)$ and the experimental order of convergence for the control as
\[e_{j}(u) = \Vert \bar u_{h_j}-\bar u_{h_{j+1}}\Vert_{L^2(\Omega)},\qquad
EOC_{j}(u) = \log_2 e_{j-1}(u) - \log_2 e_{j}(u),\]
and similarly for the state and the adjoint state. The results of the experiment are summarized in Table \ref{WWTb6.1}, where the predicted orders of convergence can be clearly noticed.

\begin{table}
  \centering
  \[
\begin{array}{cccccccc}
 j  &     e_{j}(u)  &  EOC_{j}(u) &     e_{j}(y)  &  EOC_{j}(y) &  e_{j}(\varphi)  & EOC_{j}(\varphi) \\ \hline
 3  &  1.4e-01      &      -      &  6.2e-02      &      -      &    3.6e-03       &      -           \\
 4  &  9.5e-02      &    0.5      &  1.6e-02      &    1.9      &    9.1e-04       &    2.0           \\
 5  &  5.3e-02      &    0.8      &  4.1e-03      &    2.0      &    2.3e-04       &    2.0           \\
 6  &  2.6e-02      &    1.0      &  1.0e-03      &    2.0      &    5.8e-05       &    2.0           \\
 7  &  1.3e-02      &    1.0      &  2.6e-04      &    2.0      &    1.4e-05       &    2.0           \\
 8  &  6.6e-03      &    1.0      &  6.4e-05      &    2.0      &    3.5e-06       &    2.0           \\
 9  &  3.3e-03      &    1.0      &  1.6e-05      &    2.0      &    8.9e-07       &    2.0           \\
 10 &  1.6e-03      &    1.0      &  4.0e-06      &    2.0      &    2.2e-07       &    2.0
 \end{array}
 \]
  \caption{Errors in $L^2(\Omega)$ and experimental orders of convergence}\label{WWTb6.1}
\end{table}

%

\providecommand{\bysame}{\leavevmode\hbox to3em{\hrulefill}\thinspace}
\providecommand{\MR}{\relax\ifhmode\unskip\space\fi MR }
\providecommand{\MRhref}[2]{%
  \href{http://www.ams.org/mathscinet-getitem?mr=#1}{#2}
}
\providecommand{\href}[2]{#2}

\appendix

\section{Some auxiliary results}\label{Apx}

The aim of this appendix is to prove the following theorem
\begin{theorem}\label{WWTA.9}
For every $\varepsilon > 0$ there exists $h_0>0$ such that
\[
  |[J''_h(u)-J''(u)]v^2| \leq \varepsilon v^2
\]
  for all $h<h_0$, all $u\in\uad$ and all $v\in L^2(\Omega)$.
\end{theorem}

We establish some preliminary lemmas.

\begin{lemma}\label{WWLA.4}For every $u\in\uad$ there exists a constant $C >0$ which depends in a monotonically non-decreasing way on $\Vert u\Vert_{L^2(\Omega)}$ but is independent of $h$ such that
\begin{align}
&\Vert z_{u,v}-z_{h,u,v}\Vert_{L^2(\Omega)} + h \Vert z_{u,v}-z_{h,u,v}\Vert_{H^1(\Omega)} \leq C h^2 \Vert v\Vert_{L^2(\Omega)},\label{EA.1}\\
&\Vert z_{u,v}-z_{h,u,v}\Vert_{L^\infty(\Omega)}  \leq C  h^{2-\dimension/2} \Vert v\Vert_{L^2(\Omega)}\label{EA.2}
\end{align}
for all $v\in L^2(\Omega)$, where $z_{u,v}\in H^2(\Omega)$ is the solution of \eqref{E2.5} and $z_{h,u,v}\in Y_h$ is the solution of \eqref{WWE4.1}. If $\beta < +\infty$, $C $ can be chosen independently of $u$.
\end{lemma}

\begin{proof}
  Since $z_{h,u,v}$ is not the finite element approximation of $z_{u,v}$, we introduce the intermediate function $z^h_{u,v}\in H^2(\Omega)$ solution of
\begin{align*}
& \left\{\begin{array}{l} Az^h + \displaystyle \frac{\partial a}{\partial y}(x,y_h(u))z^h + uz =- vy_h(u)\ \  \mbox{in } \Omega,\vspace{2mm}\\  \partial_{\conormal_A} z^h = 0\ \ \mbox{on }\Gamma. \end{array}\right.
\end{align*}
We have that $z_{h,u,v}$ is the finite element approximation of $z^h_{u,v}$, hence we have
\begin{align*}
&\Vert z^h_{u,v}-z_{h,u,v}\Vert_{L^2(\Omega)} + h \Vert z_{u,v}-z_{h,u,v}\Vert_{H^1(\Omega)} \leq c h^2 \Vert v\Vert_{L^2(\Omega)},\\
&\Vert z^h_{u,v}-z_{h,u,v}\Vert_{L^\infty(\Omega)}  \leq c h^{2-\dimension/2} \Vert v\Vert_{L^2(\Omega)}.
\end{align*}
On the other hand, the difference $\zeta = z_{u,v}-z^h_{u,v}$ satisfies
\begin{align*}
& \left\{\begin{array}{l} A\zeta + \displaystyle \frac{\partial a}{\partial y}(x,y_u)\zeta + u\zeta = \left(\frac{\partial a}{\partial y}(x,y_u) -\frac{\partial a}{\partial y}(x,y_h(u)) \right)z^h_{u,v}- v(y_u-y_h(u))\ \  \mbox{in } \Omega,\vspace{2mm}\\  \partial_{\conormal_A} \zeta = 0\ \ \mbox{on }\Gamma. \end{array}\right.
\end{align*}
To conclude the proof we notice that the norm in $L^\infty(\Omega)\cap H^1(\Omega)$ is bounded in terms of the norm in $L^2(\Omega)$ of the right hand side of this equation. Using assumption \eqref{A2.2}, the fact that $\Vert z^h_{u,v}\Vert_{L^\infty(\Omega)}\leq M\Vert v\Vert_{L^2(\Omega)}$ with $M$ independent of $u$ an $v$, and the error estimate for the state \eqref{WWE3.2}, we obtain the existence of a constant $C>0$ such that
\[\Vert z_{u,v}-z^h_{u,v}\Vert_{L^\infty(\Omega)\cap H^1(\Omega)}\leq C h^2.\]
and the proof is complete.
\end{proof}

The following corollary is a straightforward consequence of Lemma \ref{WWLA.4}.

\begin{corollary}\label{WWCA.5}For every $u\in\uad$ there exist a positive constants $K_1$ that depends in a monotonically non-decreasing way of $\Vert u\Vert_{L^2(\Omega)}$ but is independent of the discretization parameter $h>0$ such that
\[
    \Vert z_{h,u,v}\Vert_{L^\infty(\Omega)} +  \Vert z_{h,u,v}\Vert_{H^1(\Omega)}   \leq K_1  \Vert v\Vert_{L^2(\Omega)}
\]
  for every $h>0$ small enough. If $\beta < +\infty$, then $K_1$ is independent of $u$.
\end{corollary}

\begin{lemma}\label{WWLA.8}
For every $\varepsilon > 0$ and every $u\in\uad$ there exist $C > 0$ and $h_0>0$ such that
\[
\Vert \eta_{u,v}-\eta_{h,u,v}\Vert_{L^2(\Omega)}\leq \varepsilon \Vert v\Vert_{L^2(\Omega)}
\]
for all $0<h<h_0$ and all  $v\in L^2(\Omega)$, where $\eta_{u,v}$ is the solution of \eqref{WWE2.9} with $v_i$ replaced by $v$ and $\eta_{h,u,v}\in Y_h$ is the solution of \eqref{WWE4.6}. If $\beta < +\infty$, then $C$ and $h_0$ do not depend on $u$.
\end{lemma}
\begin{proof}Consider $\eta^h_{u,v}\in H^2(\Omega)$ the solution of
\begin{equation}\label{WWEA.5}
 \left\{\begin{array}{l} A^*\eta^h + \displaystyle \frac{\partial a}{\partial y}(x,y_h(u))\eta^h + u\eta^h\\ \hspace{1cm} \displaystyle = \left[\frac{\partial^2 L}{\partial y^2}(x,y_h(u))-\varphi_h(u) \frac{\partial^2 a}{\partial y^2}(x,y_h(u))\right] z_{h,u,v}-v\varphi_h(u)\ \  \mbox{in } \Omega,\vspace{2mm}\\  \partial_{\conormal_{A^*}} \eta^h = 0\ \ \mbox{on }\Gamma. \end{array}\right.
\end{equation}There exists a constant $K_2$ independent of $u$ and $v$ such that
\[
\Vert \eta^h\Vert_{L^\infty(\Omega)} +  \Vert \eta^h\Vert_{H^1(\Omega)}  \le  K_2 \Big\Vert \left[\frac{\partial^2 L}{\partial y^2}(x,y_h(u))-\varphi_h(u) \frac{\partial^2 a}{\partial y^2}(x,y_h(u))\right] z_{h,u,v}-v\varphi_h(u)\Big\Vert_{L^2(\Omega)}.
\]
From our assumptions and using Corollary \ref{WWCA.5} we deduce the existence of $K_3$ independent of $v$ and $h$ such that
\[
\Vert \eta^h\Vert_{L^\infty(\Omega)} +  \Vert \eta^h\Vert_{H^1(\Omega)}  \le  K_3 \Vert v\Vert_{L^2(\Omega)}.
\]
By the triangle inequality we have
\[
\Vert \eta_{u,v}-\eta_{h,u,v}\Vert_{L^2(\Omega)}\leq \Vert \eta_{u,v}-\eta^h_{u,v}\Vert_{L^2(\Omega)} + \Vert \eta^h_{u,v}-\eta_{h,u,v}\Vert_{L^2(\Omega)}
\]
Since $\eta_{h,u,v}$ is the finite element approximation of $\eta^h_{u,v}$ we have that
$\Vert \eta^h_{u,v}-\eta_{h,u,v}\Vert_{L^2(\Omega)}\leq C_1 h^2\Vert v\Vert_{L^2(\Omega)}.$
Set $h_1 = \sqrt{\frac{\varepsilon}{2C_1}}$, then
\begin{equation}\label{WWEA.6}
\Vert \eta^h_{u,v}-\eta_{h,u,v}\Vert_{L^2(\Omega)}\leq \frac{\varepsilon}{2}\Vert v\Vert_{L^2(\Omega)}\ \forall h<h_1.
\end{equation}
Let us estimate the difference $\zeta = \eta_{u,v}-\eta^h_{u,v}$. This function satisfies
\begin{align*}
  A^*\zeta  & + \displaystyle \frac{\partial a}{\partial y}(x,y_u)\zeta + u\zeta =
  \left(\displaystyle \frac{\partial a}{\partial y}(x,y_h(u))-\displaystyle \frac{\partial a}{\partial y}(x,y_u)\right)\eta^h_{u,v} \\
   & +
  \left[\frac{\partial^2 L}{\partial y^2}(x,y_u) - \frac{\partial^2 L}{\partial y^2}(x,y_h(u))\right] z_{h,u,v}
   -
  \varphi_u \left(\frac{\partial^2 a}{\partial y^2}(x,y_u)- \frac{\partial^2 a}{\partial y^2}(x,y_h(u))\right) z_{h,u,v} \\
  & -
  ( \varphi_u -\varphi_h(u)) \left(\frac{\partial^2 a}{\partial y^2}(x,y_h(u)) z_{h,u,v}+v\right) \\
   & +
  \left[\frac{\partial^2 L}{\partial y^2}(x,y_u) - \varphi_u \frac{\partial^2 a}{\partial y^2}(x,y(u))\right] (z_{u,v}-z_{h,u,v}) \mbox{ in }\Omega, \\
    \partial_{\conormal_{A^*}} \zeta & = 0\ \ \mbox{on }\Gamma.
\end{align*}
Using the Mean Value Theorem, assumptions \ref{A2.2} and \ref{WWA2.8}, the finite element error estimates for the state \eqref{WWE3.2}, the adjoint state \eqref{WWE4.8}, and the linearized state \eqref{EA.1}, we deduce the existence of $C_2>0$ such that
\begin{align*}
\Bigg\Vert  & \left(\displaystyle \frac{\partial a}{\partial y}(x,y_h(u))-\displaystyle \frac{\partial a}{\partial y}(x,y_u)\right)\eta^h_{u,v} \\
& -
( \varphi_u -\varphi_h(u)) \left(\frac{\partial^2 a}{\partial y^2}(x,y_h(u)) z_{h,u,v}+v\right) \\
& +
  \left[\frac{\partial^2 L}{\partial y^2}(x,y_u) - \varphi_u \frac{\partial^2 a}{\partial y^2}(x,y(u))\right] (z_{u,v}-z_{h,u,v})\Bigg\Vert_{L^2(\Omega)} \leq C_2 h^2 \Vert v\Vert_{L^2(\Omega)}.
\end{align*}
We put $h_2 = \sqrt{\frac{\varepsilon}{4K_2C_2}}$.

Since $\Vert y_h(u)-y_u\Vert_{L^\infty(\Omega)}\leq c h^{2-\dimension/2}$ (see Lemma \ref{WWL3.2}), using Assumptions \ref{A2.2} and \ref{WWA2.8} with $\varepsilon$ replaced by
$\displaystyle\frac{\varepsilon}{8K_1K_2}\min\{1,\frac{1}{C_\infty}\}$, where $C_\infty$ was introduced in Remark \ref{WWR2.10}, we infer the existence of $h_3>0$ such that for $h<h_3$
\begin{align*}
  &\left\Vert
  \left[\frac{\partial^2 L}{\partial y^2}(x,y_u) - \frac{\partial^2 L}{\partial y^2}(x,y_h(u))\right] z_{h,u,v}
   -
  \varphi_u \left(\frac{\partial^2 a}{\partial y^2}(x,y_u)- \frac{\partial^2 a}{\partial y^2}(x,y_h(u))\right) z_{h,u,v}
  \right\Vert_{L^2(\Omega)} \\
  &\leq \frac{\varepsilon}{8 K_1K_2} K_1\Vert v\Vert_{L^2(\Omega)} + \frac{\varepsilon}{8 K_1K_2C_\infty} K_1 C_\infty \Vert v\Vert_{L^2(\Omega)} =\frac{\varepsilon}{4K_2}  \Vert v\Vert_{L^2(\Omega)}
\end{align*}
So we have that, for $h<\min\{h_2,h_3\}$
\[\Vert \eta_{u,v}-\eta^h_{u,v}\Vert_{L^2(\Omega)} \leq K_2\left(C_2 h_2^2 + \frac{\varepsilon}{4 K_2} \right) \Vert v\Vert_{L^2(\Omega)} = \frac{\varepsilon}{2} \Vert v\Vert_{L^2(\Omega)}.\]
The result follows for $h_0 = \min\{h_1,h_2,h_3\}$ from this estimate and \eqref{WWEA.6}.
\end{proof}

\begin{proof}[Proof of Theorem \ref{WWTA.9}.]
In virtue of Lemma \ref{WWLA.8}, there exists $h_1>0$ such that $\Vert \eta_{u,v}-\eta_{h,u,v}\Vert_{L^2(\Omega)} < \frac{\varepsilon}{2 \max\{1,M_\infty\}}\Vert v\Vert_{L^2(\Omega)}$ for every $h \le h_1$, where $M_\infty$ is the constant introduced in \eqref{E2.3}. Using that $\Vert \eta_{u,v}\Vert_{L^2(\Omega)} \le K_4\Vert v\Vert_{L^2(\Omega)}$ with $K_4$ independent of $u$ and $v$, we infer
\[
\Vert \eta_{h,u,v}\Vert_{L^2(\Omega)} \le \Big(\frac{\varepsilon}{2 \max\{1,M_\infty\}} + K_4\Big)\Vert v\Vert_{L^2(\Omega)}\quad \forall h \le h_1.
\]

From the expressions obtained in \eqref{WWE2.7} and \eqref{WWE4.4}, we get with \eqref{WWE3.2}, \eqref{WWE4.8}, and Corollary \ref{WWCA.5} that
\begin{align*}
& (J''_h(u) -J''(u))v^2 =   \int_\Omega [\varphi_u z_{u,v} -\varphi_h(u) z_{h,u,v} + y_u \eta_{u,v}-y_h(u) \eta_{h,u,v}] v \dx \\
&    = \int_\Omega (\varphi_u z_{u,v} - \varphi_u z_{h,u,v})v dx + \int_\Omega (\varphi_u z_{h,u,v} - \varphi_h(u) z_{h,u,v})v\dx\\
&+ \int_\Omega (y_u \eta_{u,v}- y_u \eta_{h,u,v})v\dx + \int_\Omega (y_u \eta_{h,u,v} - y_h(u) \eta_{h,u,v}) v \dx\\
& \leq \Vert \varphi_u\Vert_{L^\infty(\Omega)} \Vert z_{u,v}-z_{h,u,v}\Vert_{L^2(\Omega)} \Vert v\Vert_{L^2(\Omega)}+
 \Vert \varphi_u - \varphi_h(u)\Vert_{L^2(\Omega)} \Vert z_{h,u,v}\Vert_{L^\infty(\Omega)} \Vert v\Vert_{L^2(\Omega)} \\
& + \Vert y_u\Vert_{L^\infty(\Omega)} \Vert \eta_{u,v}-\eta_{h,u,v}\Vert_{L^2(\Omega)} \Vert v\Vert_{L^2(\Omega)}+
     \Vert y_u-y_h(u)\Vert_{L^\infty(\Omega)} \Vert \eta_{h,u,v}\Vert_{L^2(\Omega)} \Vert v\Vert_{L^2(\Omega)}\\
&  \leq (C_\infty C  + C K_1) h^2\Vert v\Vert_{L^2(\Omega)}^2 + \Big[M_\infty \Vert \eta_{u,v}-\eta_{h,u,v}\Vert_{L^2(\Omega)} + Ch^{2-\dimension/2}\Vert \eta_{h,u,v}\Vert_{L^2(\Omega)}\Big] \Vert v\Vert_{L^2(\Omega)}\\
&\le \Big[(C_\infty C  + C K_1) h^2 + \frac{\varepsilon}{2} + Ch^{2-\dimension/2}\Big(\frac{\varepsilon}{2 \max\{1,M_\infty\}} + K_4\Big)\Big]\Vert v\Vert_{L^2(\Omega)}^2.
  \end{align*}
Taking $h_2$ such that $(C_\infty C  + C K_1) h_2^2 + Ch_2^{2-\dimension/2}\Big(\frac{\varepsilon}{2 \max\{1,M_\infty\}} + K_4\Big) = \frac{\varepsilon}{2}$,
the result follows for $h_0 = \min\{h_1,h_2\}$.
\end{proof}

\end{document}